\documentclass{amsart}
  
  \usepackage{amsthm} 
  \usepackage{amsmath} 
  \usepackage{amssymb} 
  \usepackage{mathtools} 
   \usepackage{amsfonts}
  \usepackage[dvipsnames]{xcolor}
  \usepackage[T1]{fontenc}
  
  \usepackage{wrapfig}
   
  \usepackage[scr=boondox]{mathalfa}

 \usepackage{graphicx}
  \usepackage{xspace}
  \usepackage[all]{xy}
  \usepackage{pinlabel}
  \usepackage{enumerate}
     \usepackage{overpic}
\usepackage[margin=1.3in]{geometry}
\usepackage{tikz}
\usetikzlibrary{arrows.meta,positioning,shapes.geometric}
\usepackage{float}

  \usepackage{hyperref}
\usepackage{zref-clever}
\usepackage{cleveref}

\usepackage[margin=1.3in]{geometry}
  \newcommand{\calF}{\mathcal{F}}

  \newcommand{\calL}{\mathcal{L}}
  \newcommand{\calM}{\mathcal{M}}
  \newcommand{\calN}{\mathcal{N}}
  
  \newcommand{\calP}{\mathcal{P}}

  \newcommand{\calU}{\mathcal{U}}
  \newcommand{\calV}{\mathcal{V}}

  \newcommand{\NN}{\mathbb{N}}

  \newcommand{\RR}{\mathbb{R}}

  \newcommand{\ZZ}{\mathbb{Z}}

  \newcommand{\gothic}{\mathfrak}
  \newcommand{\Ga}{{\gothic a}}
  \newcommand{\Gb}{{\gothic b}}
  \newcommand{\Gc}{{\gothic c}}
  
  \newcommand{\Ge}{{\gothic e}}

  \newtheorem{theorem}{Theorem}[section]
  \newtheorem{proposition}[theorem]{Proposition}
  \newtheorem{corollary}[theorem]{Corollary}
  \newtheorem{lemma}[theorem]{Lemma}

  \theoremstyle{definition}
  \newtheorem{definition}[theorem]{Definition}

  \newtheorem*{claim*}{Claim}

  \newtheorem*{question*}{Question}
  \newtheorem*{answer*}{Answer}
  \newtheorem*{application*}{Application}

  \theoremstyle{remark}
  \newtheorem{remark}[theorem]{Remark}
  \newtheorem*{remark*}{Remark}
  
  \DeclareMathOperator{\id}{id}

  \DeclareMathOperator{\Stab}{Stab}

  \newcommand{\Map}{\ensuremath{\operatorname{Map}}\xspace}    

  \newcommand{\param}{{\mathchoice{\mkern1mu\mbox{\raise2.2pt\hbox{$
  \centerdot$}}
  \mkern1mu}{\mkern1mu\mbox{\raise2.2pt\hbox{$\centerdot$}}\mkern1mu}{
  \mkern1.5mu\centerdot\mkern1.5mu}{\mkern1.5mu\centerdot\mkern1.5mu}}}

  \renewcommand{\setminus}{{\smallsetminus}}
  \newcommand{\ssm}{\smallsetminus}

  \newcommand{\from}{\colon\thinspace} 
  
  \newcommand{\eps}{\varepsilon}
  \newcommand{\vphi}{\varphi}
  \newcommand{\del}{\partial}

  \newcommand{\Xab}{\ensuremath{X_{\alpha,\beta}}}
  \newcommand{\Aab}{\ensuremath{A_{\alpha,\beta}}}
  
  \newcommand{\mute}[1]{}
  
\begin{document}

\title{Hyperbolicity of Big Mapping Class Groups}
\author{Assaf Bar-Natan and Yvon Verberne}
\date{\today}

\maketitle
\begin{abstract}
This paper proves that all CB-generated big mapping class groups are either coarsely bounded (and thus, trivially hyperbolic), non-trivially hyperbolic, or have infinite-coarse rank. 
\end{abstract}

\vspace{-0.25cm}

\section{Introduction}

Let $\Sigma$ be a connected, orientable, infinite-type surface such that the end space of $\Sigma$ is \textit{stable}
(see Section~\ref{Subsec:CommonProperties}). Let $\Map(\Sigma)$ denote the mapping class group of $\Sigma$, defined to be the group of orientation preserving homeomorphisms up to isotopy.
Further, suppose that $\Map(\Sigma)$ admits a \textit{CB-generating set} (see Section~\ref{Subsec:CommonProperties}).

This preprint, together with the work of Mann-Rafi \cite{MR} and Bar-Natan--Qing--Rafi \cite{BNQR} proves the following main result:

\begin{theorem}\label{MainThm}
Suppose that $\Sigma$ is stable and let $\Map(\Sigma)$ have a coarsely bounded generating set. Then $\Map(\Sigma)$ has infinite coarse rank or is hyperbolic. Moreover, If $\Map(\Sigma)$ is hyperbolic and is not coarsely bounded, then $\Sigma$ is either a broom, brush, or tripod surface.
\end{theorem}

We summarize the coarse-rank landscape in the following flowchart:
\providecommand{\flowchartscale}{0.7}

\begin{figure}[htbp]
\centering
\scalebox{\flowchartscale}{%
\begin{tikzpicture}[
    font=\small,
    >={Stealth[round]},
    base/.style   ={draw, align=center, inner sep=6pt, line width=.7pt,
                    rounded corners=3pt,
                    execute at begin node={\hyphenpenalty=10000\relax}},
    setup/.style  ={base, dashed, draw=black!45, fill=black!3,
                    text width=76mm, font=\small\itshape},
    decision/.style={base, fill=blue!8,   draw=blue!55!black,   text width=42mm},
    rank/.style   ={base, fill=red!8,     draw=red!55!black,    text width=42mm},
    cb/.style     ={base, fill=green!12,  draw=green!45!black,  text width=34mm},
    hyp/.style    ={base, fill=orange!14, draw=orange!70!black, text width=52mm},
    arrow/.style  ={->, line width=.7pt, draw=black!70},
    lbl/.style    ={font=\scriptsize, inner sep=2pt,
                    fill=white, fill opacity=.85, text opacity=1},
  ]

  \node[setup]    (as) at (0,1.9)
        {$\Sigma$ is stable and
         $\Map(\Sigma)$ is CB-generated.};

  \node[decision] (q1) at (0,0)
        {Does $\Sigma$ admit essential shifts?};

  \node[rank]     (t1) at (6.2,0)
        {$\Map(\Sigma)$ has infinite coarse rank\\[2pt]
         {\scriptsize \cite{GRV}}};

  \node[decision] (q2) at (0,-2.7)
        {What is $\zeta(\Sigma)$?};

  \node[rank]     (t2) at (6.2,-2.7)
        {$\Map(\Sigma)$ has infinite coarse rank\\[2pt]
         {\scriptsize \cite{BNQR}}};

  \node[cb]       (t3) at (-6.2,-2.7)
        {$\Map(\Sigma)$ is CB\\[2pt]
         {\scriptsize \cite{MR}, \cite{BNQR}}};

  \node[decision] (qx) at (0,-6.0)
        {What is $\xi(\Sigma)$?};

  \node[rank]     (t4) at (6.2,-6.0)
        {$\Map(\Sigma)$ has infinite coarse rank\\[2pt]
         {\scriptsize arbitrarily many disjoint\\ non-peripheral curves}};

  \node[decision] (q3) at (-6.2,-6.0)
        {Does $\Sigma$ have finite genus, or some maximal ends
         accumulated by genus and some not?};

  \node[rank]     (t6) at (-6.2,-9.6)
        {Arbitrarily many\\ non-peripheral curves};

  \node[hyp]      (t5) at (0,-9.6)
        {$\Sigma$ is a broom / brush / tripod\\[2pt]
         and $\Map(\Sigma)$ is hyperbolic};

  \draw[arrow] (as) -- (q1);

  \draw[arrow] (q1) -- node[lbl,above]{yes} (t1);
  \draw[arrow] (q1) -- node[lbl,right]{no}  (q2);

  \draw[arrow] (q2) -- node[lbl,above]{$\zeta(\Sigma)>3$}   (t2);
  \draw[arrow] (q2) -- node[lbl,above]{$\zeta(\Sigma)\le 2$}(t3);
  \draw[arrow] (q2) -- node[lbl,right,align=left]{$\zeta(\Sigma)=3$} (qx);

  \draw[arrow] (qx) -- node[lbl,above]{$\xi(\Sigma)>3$} (t4);
  \draw[arrow] (qx) -- node[lbl,above]{$\xi(\Sigma)=3$} (q3);

  \draw[arrow] (q3) -- node[lbl,right]{yes} (t6);
  \draw[arrow] (q3.south east) -- node[lbl,above right]{no} (t5.north west);

\end{tikzpicture}%
}
\caption{Coarse geometry of $\Map(\Sigma)$ for $\Sigma$ stable with
         $\Map(\Sigma)$ CB-generated.}
\label{fig:coarse-geometry-map}
\end{figure}
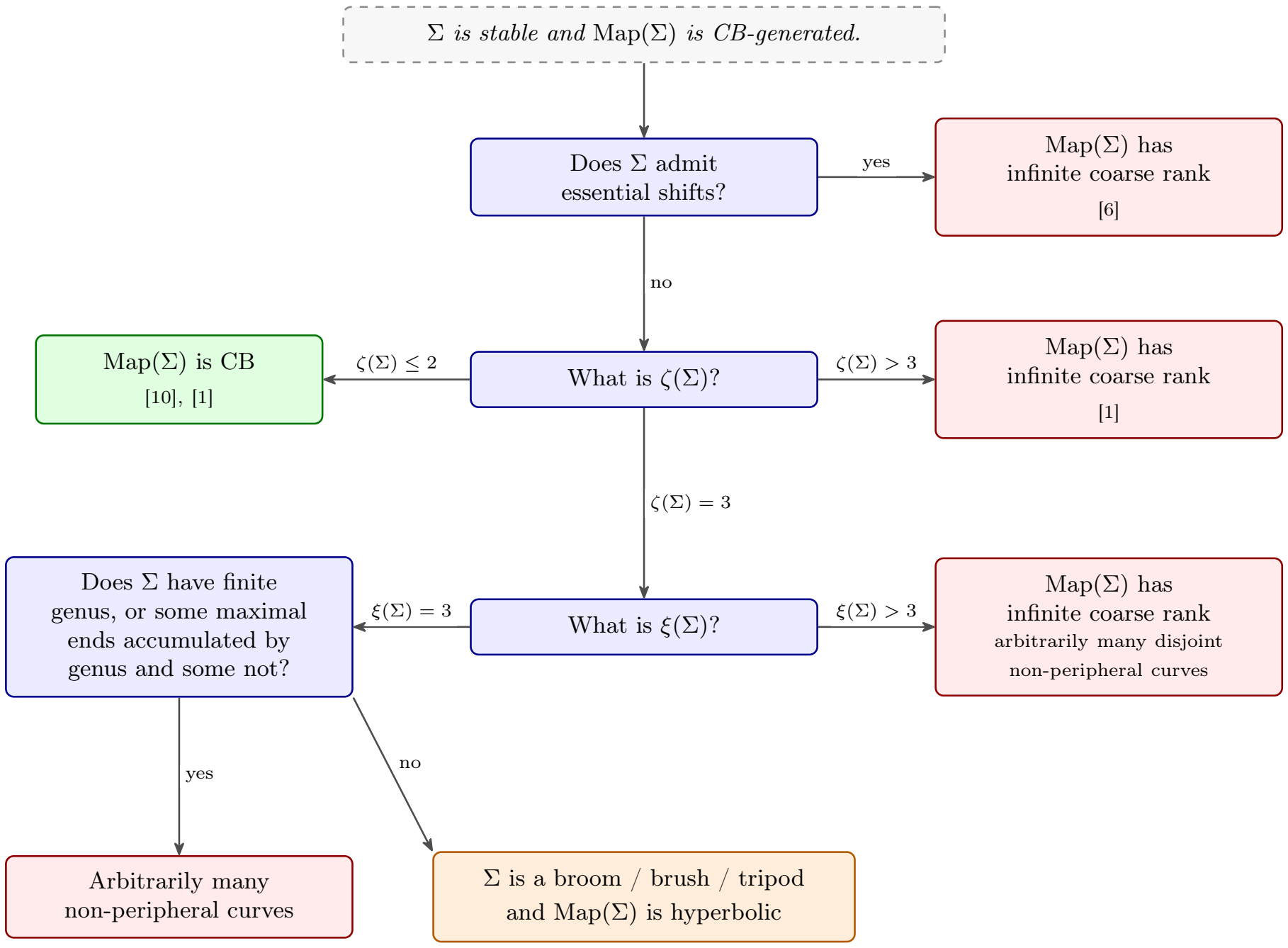

The question of which big mapping class groups which have a coarsely bounded generated set, otherwise referred to as \textit{CB-generated} big mapping class groups, are hyperbolic was first posed to the authors by Horbez--Qing--Rafi in 2021, but later appears in \cite{BKR} as Problem 2.46(a) posed by Abbott--Mann, and within Conjecture B in Bar-Natan--Qing--Rafi \cite{BNQR}. 
This result not only answers which big mapping class groups are hyperbolic, but also completes the classification of CB-generated big mapping class groups with infinite coarse rank introduced by Bar-Natan--Qing--Rafi \cite{BNQR}.

\subsection{Prior results}
With respect to the hyperbolicity of big mapping class groups,
the work of Mann--Rafi classified for which surfaces $\Sigma$, $\Map(\Sigma)$ is \textit{coarsely bounded}, abbreviated CB, (see Section~\ref{Subsec:CommonProperties}).
Any big mapping class group which is CB is quasi-isometric to a point, implying that the group itself is \textit{trivially hyperbolic}. 
If a big mapping class group is not CB but is hyperbolic, we will say that the group is \textit{non-trivially hyperbolic}.
The first example of a big mapping class group which is non-trivially hyperbolic was provided by Schaffer-Cohen, who showed that the mapping class group of the plane minus a Cantor set is hyperbolic \cite{SC}.
Schaffer-Cohens approach was to prove that the mapping class group of a plane minus a Cantor set is quasi-isometric to the loop graph, a graph which was proven by Bavard to be hyperbolic \cite{Bav}.
The approach of Schaffer-Cohen is modified in this work to prove that the mapping class group for three further classes of surfaces is non-trivially hyperbolic.

On the other end of the spectrum, there has been a considerable amount of progress made with regards to proving that big mapping class groups fail to be hyperbolic.
One direction is to prove the existence of two disjoint nondisplaceable subsurfaces. A subsurface $S \subset \Sigma$ is said to be nondisplaceable if for any $f \in \Map(\Sigma)$, $f(S) \cap S \neq \emptyset$. See Section 2 of \cite{MR} for details.

Another related direction is to determine whether the big mapping class group has infinite asymptotic dimension. We note that infinite asymptotic dimension isn't a sufficient condition to prove that the mapping class group fails to be hyperbolic. 
The work of Grant--Rafi--Verberne proves that the existence of an \textit{essential shift} (see Section~\ref{Subsec:CommonProperties}) implies the existence of high dimensional quasi-flats, and thus, that the mapping class group has infinite asymptotic dimension \cite{GRV}. 
When the mapping class group contains an essential shift, we remark that one can conclude that the mapping class group fails to be hyperbolic.
In work of Kopreski--Shaji, arc-and-curve models for locally bounded Polish subgroups of big mapping class groups are constructed, and they show that the asymptotic dimension of a CB-generated big mapping class group is infinite unless the big mapping class group is CB \cite{KS}.
In this case, one is unable to conclude whether the big mapping class group is hyperbolic.

On \cite{BNQR}, Bar-Natan--Qing--Rafi, study the \textit{coarse rank} big mapping class groups, where coarse rank is defined to be the largest $n$ such that $\ZZ^n$ quasi-isometrically embeds into the group.
An example of mapping class groups with infinite coarse ranks are those which contain an essential shift. 
Additionally, the authors introduce a numerical invariant $\zeta$ which measures the end complexity of a surface with CB-generated big mapping class group. See Section~\ref{Subsec:CommonProperties} for a formal definition.
Within this setting, they are able to provide a full classification of mapping class groups which have infinite coarse rank, except in the case where $\zeta(\Sigma) = 3$.
This paper studies the surfaces for which this numerical constant equals 3. To aid with the classification, we introduce a 
modified invariant, $\xi$ which refines the invariant $\zeta$.

\subsection{Overview}

The surfaces we study in this paper fall into six categories: the broom, brush and tripod surfaces, and the closely related broom-like, brush-like and tripod-like surfaces. See Section~\ref{subsec:broombrushtripod} for full definitions of these six categories of surfaces. 
The authors would like to remark that the definition of a tripod surface in this paper is a more broad than the standard definition of a surface with no punctures and exactly three ends accumulated by genus.

In the case of the broom, brush and tripod surfaces, we construct a graph where the vertices consist of nondisplaceable pairs of pants, and edges correspond to the underlying pairs of pants having exactly two points of intersection.
We show that this graph is hyperbolic, and in the case where $\Sigma$ is either the broom, brush or tripod surface, then the mapping class group is quasi-isometric to the nondisplaceable pants graph.

In the case of the broom-like, brush-like and tripod-like surfaces, we show that for any $k \in \NN$, there exists $k$ pairwise disjoint non-peripheral curves in each surface.
Work of Bar-Natan--Qing--Rafi proves that this is a sufficient condition to conclude that the mapping class group of each of these surfaces has infinite coarse rank.

\subsection{Outline of paper}
In Section 2, we provide the readers with the full definitions of the surfaces we will study in this paper, including a description of the endspace for three classes of surfaces.
In Section 3, we introduce the nondisplaceable pants graph, and prove properties which hold for vertices which are adjacent in the graph.
In Section 4, we introduce an arc and loop graph, which we show is hyperbolic using the Guessing Geodesics Lemma.
In Section 5, we use the arc and loop graph together with the Guessing Geodesics Lemma to prove that the nondisplaceable pants graph is non-trivially hyperbolic for the broom, brush and tripod surfaces.
In Section 6, we prove that there exists a quasi-isometry between the mapping class group of the broom, brush and tripod surfaces and the nondisplaceable pants graphs for these surfaces.
In Section 7, we prove that surfaces with $\xi > 3$ and that the broom-like, brush-like and tripod-like surfaces all have infinite coarse rank, which completes the proof of the main theorem.\\

\noindent \textbf{Acknowledgements.}
The authors are grateful to Anschel Schaffer-Cohen for preliminary discussions at the inception of this project, for discussions related to the nondisplaceable pants graph and the arc and loop graph, for providing the proof of Proposition~\ref{prop:qi} for the broom surface, and for providing the authors with the tools to prove the result for both the brush and tripod surfaces.
The authors are also grateful to Camille Horbez, Yulan Qing, and Kasra Rafi for both suggesting this problem, and for sharing their preliminary notes pertaining to which CB-generated big mapping class groups are hyperbolic.
The authors would also like to thank Kasra Rafi for conversations surrounding the hyperbolicity of the arc and loop graph, specifically Lemma~\ref{lem:GuessingGeodesicslocalcentralpants}.
The authors would further like to thank Carolyn Abbott, Jesús Hernández Hernández, Camille Horbez, Dan Margalit, Yulan Qing, Lorenzo Ruffoni, and Ferrán Valdez for helpful discussions.
Finally, the authors would like to thank the organizers of the conferences ``Big Mapping Class Groups and Diffeomorphism Groups" at CIRM and ``Geometric Topology in New York" at Columbia University for both their hospitality and providing an opportunity for the authors to collaborate in person. Y.V. was  partially supported by an NSERC-PDF Fellowship and by an NSERC Discovery grant RGPIN 05587.\\

\noindent \textbf{AI Disclosure.}
AI was used to proofread a draft of the manuscript.
Kasra Rafi sent an AI generated report on an initial draft of the paper, which the authors responded to.

\section{Broom-Like, Brush-Like, and Tripod-Like Surfaces}
In this section, we wish to discuss the definitions and topological properties for the broom-like, brush-like, and tripod-like surfaces which were introduced in the the introduction. 
The first subsection will formally define each of the three surfaces that we will prove have hyperbolic big mapping class group, then discuss some topological properties of these surfaces in the next subsection.
The final subsection will give a brief discussion of the remaining broom-like, brush-like, and tripod-like surfaces which we will show have infinite coarse rank.

\subsection{Four Common Properties}\label{Subsec:CommonProperties}

Before we provide the reader with the formal definitions of our surfaces of interest, we set the stage by defining some of the important properties each of our surfaces share. We will consider surfaces that

\begin{enumerate}
\item Are stable;
\item Have CB-generated mapping class groups;
\item Admit no essential shifts; and
\item Have end complexity at least $3$.
\end{enumerate}

To begin, a subset $X$ of a group $G$ is \textit{coarsely bounded}, abbreviated
CB, if every compatible left-invariant metric on $G$ gives $X$ finite diameter.
We see that a group is \textit{CB-generated} if it admits a CB-generating set.

The \textit{space of ends} of a surface $\Sigma$, denoted $E = E(\Sigma)$, is
defined as $E(\Sigma) = \varprojlim_{K\subset \Sigma} \pi_0(\Sigma\ssm K)$,
where $K$ ranges over all compact subsets of $\Sigma$, and if $K_i\subset K_j$,
then the map in the projective limit is the inclusion
$\vphi_{ij}:\pi_0(\Sigma\ssm K_j)\hookrightarrow \pi_0(\Sigma\ssm K_i)$.
Intuitively, each end corresponds to a way one can leave every compact
subsurface of $\Sigma$.  For $x \in End(\Sigma)$, we call a neighborhood $U$ of
$x$ \textit{stable} if for any smaller neighbourhood $U' \subset U$ of $x$,
there is a homeomorphism from $U$ to $U'$ fixing $x$.

A \textit{shift map} is a homeomorphism which translates the ends of the
surface along a predefined infinite strip.  More precisely, consider an
infinite strip in $\RR^2$ with $\ZZ$ acting on the strip by translations. Cut
out an equivariant family of disks and attach identical surfaces (possibly of
infinite type) to the boundaries of these disks.  Notice that $\ZZ$ still acts
by translations on the strip.  Embedding this into a larger strip $\sigma$, we
can construct a homeomorphism $h_{\sigma}$ of $\sigma$ that acts as described
above in the smaller strip, but the restriction of $h_{\sigma}$ to $\partial
\sigma$ is the identity. 
A copy of $\sigma$ where the ends of $\sigma$ exit the different ends of
$\Sigma$ is a shift map.  We say a shift map $h_{\sigma}$ is \textit{essential}
if $\langle h_{\sigma}\rangle$, is not a coarsely bounded subgroup of
$\Map(\Sigma)$.  Equivalently, an essential shift is a countable shift map that
cannot be pushed into the neighborhood of an end.

Mann--Rafi define a partial ordering on the space of ends by setting
$x \prec y$, for $x, y \in E(\Sigma)$, if every clopen neighbourhood of $y$ contains a clopen set
homeomorphic to a neighbourhood of $x$ \cite{MR}. Two ends are equivalent, or are of the
same type if $x \preceq y$ and $y \preceq x$. This partial order has maximal elements, and we let
$\calM(\Sigma)$ denote the space of \textit{maximal ends} of $\Sigma$. 
In Bar-Natan--Qing--Rafi, the authors introduce the \textit{end complexity} of a surface $\Sigma$ which is defined as follows:
    \begin{equation*}
        \begin{aligned}
            \zeta(\Sigma) &= 2 \cdot \# (E(x) \text{ where } E(x) \subset \calM(End(\Sigma)) \text{ is a Cantor set}) +\\
            &+ \#(\text{of isolated points in } \calM(End(\Sigma)))+\\
            &+ \#(\text{isolated maximal ends that are unique maximal accumulation points}).
        \end{aligned}
    \end{equation*}
The end complexity is a numerical invariant of a surface $\Sigma$.

\subsection{Broom, brush, and tripod surfaces}\label{subsec:broombrushtripod}

In this subsection, we define the broom, brush and tripod surfaces.

\begin{definition}[Broom Surface]\label{def:broom}
    A surface $\Sigma$ is a \emph{broom surface} if:
    \begin{itemize}
        \item $\Sigma$ is stable and $\Map(\Sigma)$ is CB-generated;
				\item $\mathcal{M}(\Sigma)$ consists of a disjoint union of a Cantor
					set and a singleton;
				\item $\Sigma$ admits no essential shifts;
				\item Every end that is not maximal accumulates to every maximal end;
					and
				\item $\Sigma$ has zero or infinite genus; in the latter case, every
					maximal end is accumulated by genus
    \end{itemize}
\end{definition}

\begin{definition}[Brush Surface]\label{def:brushsurface}
    A surface $\Sigma$ is a \emph{brush surface} if:
    \begin{itemize}
				\item $\Sigma$ is stable and $\Map(\Sigma)$ is CB-generated;
				\item $\mathcal{M}(\Sigma)$ consists of two elements, $x$ and $y$;
				\item $\Sigma$ admits no essential shifts;
				\item Every equivalence class of ends that accumulates to $y$ accumulates to $x$;
				\item There exists precisely one equivalence class of ends that 
					accumulates to $x$ but does not accumulate to $y$; and
				\item $\Sigma$ has zero or infinite genus; in the latter case, every
					maximal end is accumulated by genus
    \end{itemize}
\end{definition}

\begin{definition}[Tripod Surface]\label{def:tripod}
    A surface $\Sigma$ is a \emph{tripod} if:
    \begin{itemize}
        \item $\Sigma$ is stable and $\Map(\Sigma)$ is CB-generated;
        \item The surface $\Sigma$ has exactly $3$ maximal ends, all equivalent;
        \item $\Sigma$ admits no essential shifts; and
        \item $\Sigma$ has zero or infinite genus; in the latter case, every maximal end is accumulated by genus.
    \end{itemize}
\end{definition}

\begin{definition}

A surface is called \textit{broom-like} (respectively, \textit{brush-like}, \textit{tripod-like}) if it satisfies all the requirements of being a broom (respectively, brush, tripod) with the exception of the condition on genus. In other words, a surface is \textit{broom-like} (respectively, \textit{brush-like}, \textit{tripod-like}) if it satisfies all the requirements of being a broom (respectively, brush, tripod) but has finite genus, or it satisfies all the requirements but
if it has at least one maximal end accumulated by genus and at least one that is not accumulated by genus.
\end{definition}

\subsection{Endspace structure of the broom, brush and tripod surfaces}

In this subsection, we describe the endspace structure of the broom, brush and tripod surfaces, as well as set the notation we will use when working with these surfaces.

\subsubsection{Endspace structure of broom surfaces}\label{subsub:endspacebroom}

Note first that by Definition~\ref{def:broom}, $\calM(\Sigma) = \{\Ga\} \sqcup X$, where $X$ is homeomorphic to a Cantor set.
For any end $e$ of $\Sigma$, if $e$ is not equivalent to one of the maximal ends then it is a predecessor of all the maximal ends. In particular, by tameness there is a finite set $\{e_i\}_{i=1}^n$ of immediate predecessors to the maximal ends such that every immediate predecessor is equivalent to some $e_i$.

By tameness, each $e_i$ has a stable neighborhood $U_i$. Because the equivalence class $E(e_i)$ is not maximal, it cannot be finite. If $E(e_i)$ were countable, it would induce an essential shift, which is ruled out by Definition~\ref{def:broom}. Thus $E(e_i)$ must be uncountable, and $U_i \cap E(e_i)$ is a Cantor set. Then $\calU = \bigcup_{i=1}^n U_i$ is a clopen set of ends containing a copy of every non-maximal end. What's more, if $V$ is any clopen set of ends of $\Sigma$ that does not contain a maximal end, and $\calU'$ is a clopen set of ends homeomorphic to $\calU$, $V \cup \calU'$ is itself homeomorphic to $\calU$.\\

\subsubsection{Endspace structure of brush surfaces}\label{subsub:endspacebrush}

Note first that by Definition~\ref{def:brushsurface}, $\calM(\Sigma) = \{\Ga, \Gb\}$, and we will denote by $E_{\Gc}$ the equivalence class of ends that is an immediate predecessor of $\Ga$ but not $\Gb$.
Note first that by Definition~\ref{def:brushsurface}, for any end $e$ of $\Sigma$, if $e$ is not equivalent to one of the maximal ends then it is a predecessor of all the maximal ends or it is an immediate predecessor of $\Ga$ but does not accumulate to $\Gb$. In the former case, by tameness there is a finite set $\{e_i\}_{i=1}^n$ of immediate predecessors to the maximal ends such that every immediate predecessor is equivalent to some $e_i$. 

Again by tameness, each $e_i$ has a stable neighborhood $U_i$. 
Because the equivalence class $E(e_i)$ is not maximal, it cannot be finite. If $E(e_i)$ were countable, it would induce an essential shift, which is ruled out by Definition~\ref{def:brushsurface}. 
Thus $E(e_i)$ must be uncountable, and $U_i \cap E(e_i)$ is a Cantor set. 
Then $\calU = \bigcup_{i=1}^n U_i$ is a clopen set of ends containing a copy of every non-maximal end which is an immediate predecessor of both $\Ga$ and $\Gb$. 
What's more, if $W$ is any clopen set of ends of $\Sigma$ that does not contain a maximal end but is an immediate predecessor of both $\Ga$ and $\Gb$, and $\calU'$ is a clopen set of ends homeomorphic to $\calU$, $W \cup \calU'$ is itself homeomorphic to $\calU$.

In the latter case, if an end $e$ of $\Sigma$ is not a maximal end and it is an immediate predecessor of $\Ga$ but does not accumulate to $\Gb$, then by tameness there is a finite set $\{e_i\}_{i=1}^n$ of immediate predecessors to the maximal end $\Ga$ such that every immediate predecessor is equivalent to some $e_i$. 
Again by tameness, each $e_i$ has a stable neighborhood $V_i$. Because the equivalence class $E(e_i)$ is not maximal, it cannot be finite. 
If $E(e_i)$ were countable, it would induce an essential shift, which is ruled out by Definition~\ref{def:brushsurface}. 
Thus $E(e_i)$ must be uncountable, and $V_i \cap E(\Ge)$ is a Cantor set. Then $\calV = \bigcup_{i=1}^n V_i$ is a clopen set of ends containing a copy of every non-maximal end which are immediate predecessors to the maximal end $\Ga$. 
What's more, if $W$ is any clopen set of ends of $\Sigma$ that does not contain a maximal end but is an immediate predecessor of $\Ga$ but does not accumulate to $\Gb$, 
and $\calV'$ is a clopen set of ends homeomorphic to $\calV$, $W \cup \calV'$ is itself homeomorphic to $\calV$.

\subsubsection{Endspace structure of tripod surfaces}\label{subsub:endspacetripod}
Note first that by Definition~\ref{def:tripod}, $\calM(\Sigma) = \{ \Ga, \Gb, \Gc\}$ all of which are equivalent.
Let $e$ be an end of $\Sigma$ which is not equivalent to one of the maximal ends. 
Then $e$ must be a predecessor of each of the maximal ends of the tripod. 
By tameness, there is a finite set $\{e_i\}_{i=1}^n$ of immediate predecessors to the maximal ends such that every immediate predecessor is equivalent to some $e_i$. Again by tameness, each $e_i$ has a stable neighborhood $U_i$. Because the equivalence class $E(e_i)$ is not maximal, it cannot be finite. If $E(e_i)$ were countable, it would induce an essential shift, which is ruled out by Definition~\ref{def:tripod}. Thus $E(e_i)$ must be uncountable, and $U_i \cap E(e_i)$ is a Cantor set. Then $\calU = \bigcup_{i=1}^n U_i$ is a clopen set of ends containing a copy of every non-maximal end. What's more, if $W$ is any clopen set of ends of $\Sigma$ that does not contain a maximal end, and $\calU'$ is a clopen set of ends homeomorphic to $\calU$, $W \cup \calU'$ is itself homeomorphic to $\calU$.

\section{Nondisplaceable pants graph}

In this section, we introduce the nondisplaceable pants graph.
Additionally, we will introduce the language used to describe the different boundary components of a nondisplaceable pair of pants.
In later sections, this graph will be shown to be $\delta$-hyperbolic and quasi-isometric to the broom, brush and tripod surfaces.

\subsection{Boundary components of a nondisplaceable pair of pants}\label{Subsec:BoundaryComponentsNDPants}

Throughout this subsection, let $\Sigma$ be either a broom, a brush, or
a tripod surface.  When $\Sigma$ is a broom, we denote $\calM(\Sigma) = \{\Ga\}
\sqcup X$ where $X$ is homeomorphic to a Cantor set, and let $\Gb$ be a prescribed element of X. 
When $\Sigma$ is a brush, we denote $\calM(\Sigma) = \{\Ga,\Gb\}$, and denote by $E_{\Gc}$ the equivalence class of
ends that converges to $\Ga$ but not $\Gb$. When $\Sigma$ is
a tripod, we denote $\calM(\Sigma) = \{\Ga,\Gb, \Gc\}$.

To begin our analysis of the boundary components of a nondisplaceable pair of pants, denoted $P$, for the broom, brush and tripod surfaces, by introducing notation which allows us differentiate the three boundary components of $P$.

\begin{definition}
Let $c \subset \Sigma$ be a simple closed curve.
\begin{itemize}
\item If $\Sigma$ is a broom we say that $c$ is an \textit{$\alpha$-curve} if
	$c$ separates $\Ga$ from $X$. We say that $c$ is a \textit{$\beta$-curve} if
		both components of $\Sigma \setminus c$ contain a subset of $X$ and $c$
		separates $\Ga$ from $\Gb$. Suppose that $\alpha$ denotes an $\alpha$-curve
		and $\beta$ denotes a $\beta$-curve. Then a $\gamma$-curve is any curve
		which forms a nondisplaceable pair of pants with $\alpha$ and $\beta$.
\item If $\Sigma$ is a brush, we say that $c$ is an \textit{$\alpha$-curve} if
	$c$ separates $\Gb$ from $\{\Ga\} \cup E_{\Gc}$. We say that $c$ is
		a \textit{$\beta$-curve} if both components of $\Sigma \setminus c$ contain
		a subset of $E_{\Gc}$, and $c$ separates $\Ga$ from $\Gb$. Suppose that
		$\alpha$ denotes an $\alpha$-curve and $\beta$ denotes a $\beta$-curve.
		Then a $\gamma$-curve is a curve which forms a nondisplaceable pair of
		pants with $\alpha$ and $\beta$.
\item If $\Sigma$ is a tripod, we say that $c$ is an \textit{$\alpha$-curve} if
	$c$ separates $\Ga$ from $\{\Gb,\Gc\}$. We say that $c$ is
		a \textit{$\beta$-curve} if $c$ separates $\Gb$ from $\{\Ga,\Gc\}$. We say
		that $c$ is a \textit{$\gamma$-curve} if $c$ separates $\Gc$ from $\{\Ga,
		\Gb\}$.
\end{itemize}
\end{definition}

\noindent We note to the reader that by definition, a $\gamma$ curve is unable to be of either $\alpha$ or $\beta$ type.

As a final piece of notation, suppose that $P = \{a,b,c\}$ is a pair of pants in a surface $\Sigma$ and let $S_a$ denote the subsurface of $\Sigma \ssm a$ which does not contain either the boundary components $b$ or $c$.
We may similarly define $S_b$ and $S_c$.

The following is an immediate consequence of the classification of noncompact surfaces by Richards \cite{R}.

\begin{lemma}\label{lem:changeofpants}
    Suppose that $K_1, K_2$ are subsurfaces of $\Sigma$ homeomorphic to the same surface $K$. Then 
    there exists $f \in Homeo(\Sigma)$ such that 
    $f(K_1) = K_2$ if and only if the complementary 
    components of $K_1$ are all homeomorphic to the 
    complementary components of $K_2$.
\end{lemma}

We begin by showing that if $\Sigma$ is a broom, brush, or tripod surface, then a pair of pants is nondisplaceable if and only if it contains both an $\alpha$-curve and a $\beta$-curve.

\begin{lemma}\label{lem:nondisplaceableimpliesalphabeta}
    Let $\Sigma$ be a broom, brush, or tripod  surface, and $P = \{a,b,c\}$ a pair of pants in $\Sigma$. 
    Then $P$ is nondisplaceable if and only if it contains an $\alpha$-curve and a $\beta$-curve.
\end{lemma}

\begin{proof}
To begin, we remark that if $P$ is nondisplaceable, then every cuff of $P$ must be separating. 
Suppose not.
Notice first, that since we are only considering the broom, brush and tripod surfaces that the existence of a non-separating curve implies the existence of infinite genus.
Indeed, on planar surfaces, every curve is separating.
Now notice that if two cuffs of $P$ are separating in $\Sigma$, then the third cuff must also be separating since $\partial P = a \cup b \cup c$ bounds $P$.
Thus, if not all cuffs of $P$ are separating, then at least two cuffs of $P$ are non-separating, which implies that $\Sigma \setminus P$ has at most two connected components.
Since every maximal end of $\Sigma$ is accumulated by genus, at least one component $U$ of $\Sigma \setminus P$ has infinite genus, allowing us to find a disjoint pair of pants $P' \subset U$ whose complement has the same end-space and genus partition as $\Sigma \setminus P$.
By Lemma~\ref{lem:changeofpants}, $P$ is displaceable, and thus we have that every cuff of $P$ is separating.

We will now prove the statement for the broom surface.
We begin with the forward direction and suppose that $P$ is a non-displaceable subsurface.
We further suppose that $P$ does not contain an $\alpha$-curve.
In the case of the broom surface, this implies that no boundary curve of $P$ separates the singleton maximal end from the maximal ends homeomorphic to a Cantor set.
In particular, each boundary curve of $P$ contains some subset of the maximal ends which are homeomorphic to a Cantor set.
Let $a$ denote the boundary curve of $P$ containing the singleton maximal end, which also includes some subset of the maximal ends which are homeomorphic to a Cantor set.
By the classification of noncompact surfaces by Richards, there exists a pair of pants with boundary components $P' = \{a', b', c'\}$ which lies strictly within $S_a$ and cuts off complementary components such that $S_{a'} \cong S_a$, $S_{b'} \cong S_b$, and $S_{c'} \cong S_c$.
This implies by Lemma~\ref{lem:changeofpants} that there exists a mapping class $f \in \Map(\Sigma)$ such that $f(P) = P'$ and thus $f(P) \cap P = \emptyset$, a contradiction.

We now suppose that $P$ contains an $\alpha$-curve but does not contain a $\beta$-curve.
This implies that $P$ has a boundary curve, denoted $a$, which separates the singleton maximal end from the subset of maximal ends homeomorphic to a Cantor set, a curve $b$ which contains the entire subset of maximal ends homeomorphic to a Cantor set, and a curve $c$ which contains no maximal ends.
Recall that in the case of the broom, an end which is not maximal must accumulate to every maximal end.
Therefore, in the subsurface $S_a$, there must exist a neighborhood equivalent to $S_c$ since the ends in $c$ accumulate to the singleton maximal end.
By the classification of noncompact surfaces by Richards, there exists a pair of pants with boundary components $P' = \{a', b', c'\}$ which lies strictly within $S_a$ and cuts off complementary components such that $S_{a'} \cong S_a$, $S_{b'} \cong S_b$, and $S_{c'} \cong S_c$.
This implies by Lemma~\ref{lem:changeofpants} that there exists a mapping class $f \in \Map(\Sigma)$ such that $f(P) = P'$ and thus $f(P) \cap P = \emptyset$, a contradiction.

For the other direction, we suppose that $P = \{a,b,c\}$ contains both an $\alpha$-curve and a $\beta$-curve. 
We let the $\alpha$-curve be denoted by $a$ and the $\beta$-curve be denoted by $b$.
By way of contradiction, we suppose that $P$ is displaceable.
Since $P$ is displaceable, this means that there exists some $f \in \Map(\Sigma)$ such that $f(P) \cap P = \emptyset$.

If $f(P) \subset S_{a}$, then since $a$ was an $\alpha$-curve $f(P)$ is unable to contain a $\beta$-curve.
Indeed, in the case of the broom surface, since $a$ separated the singleton maximal end from the subset of maximal ends homeomorphic to a Cantor set, the entire subset of ends homeomorphic to a Cantor set lies to one side of a boundary component of $f(P)$.

If $f(P) \subset S_{b}$, then since $b$ was an $\beta$-curve $f(P)$ is unable to contain a $\alpha$-curve.
Indeed, in the case of the broom surface, since $b$ contains a subset of the maximal ends homeomorphic to a Cantor set, no boundary component of $f(P)$ is able to separate the singleton maximal end from the maximal ends homeomorphic to a Cantor set.

If $f(P) \subset S_{c}$, then since $b$ was neither an alpha curve nor a $\beta$-curve $f(P)$ is unable to an $\alpha$-curve.
Indeed, in the case of the broom surface, since $c$ contains a subset of the cantor set of ends, no boundary component of $f(P)$ is able to separate the singleton maximal end from the maximal ends homeomorphic to a Cantor set, completing the proof for the broom surface.

The proof above follows similarly for the brush surface, and we leave this as an exercise for the reader.

For the tripod surface, begin with the forward direction and suppose that $P = \{a,b,c\}$ is a nondisplacable pair of pants.
Since the three cuffs of $P$ are disjoint separating curves, $\Sigma \ssm P$ consists of $P$, and three pairwise disjoint subsurfaces $S_a, S_b$ and $S_c$.
By way of contradiction, suppose that it is not the case that each of $S_a, S_b$ and $S_c$ contains a single maximal end.
This implies, without loss of generality, that $S_a$ contains no maximal ends, and at least one, say $S_b$ contains at least two maximal ends by the Pigeonhole Principle.
By the classification of noncompact surfaces by Richards, this implies that there exists a pair of pants with boundary components $P' = \{a', b', c'\}$ which lies strictly within $S_b$ and cuts off complementary components such that $S_{a'} \cong S_a$, $S_{b'} \cong S_b$, and $S_{c'} \cong S_c$.
This implies by Lemma~\ref{lem:changeofpants} that there exists a mapping class $f \in \Map(\Sigma)$ such that $f(P) = P'$ and thus $f(P) \cap P = \emptyset$, a contradiction.

For the other direction, suppose that $P = \{a,b,c\}$ contains both an alpha curve and a $\beta$ curve.
If $f(P) \subset S_a$, then by definition notice that the ends $\Gb$ and $\Gc$ lie in the same complementary component of $\Sigma \ssm f(P)$, implying that $P \not \cong f(P)$.
The same argument holds if $f(P) \subset S_b$ or if $f(P) \subset S_c$, completing the proof of the claim.
\end{proof}

We now define the nondisplaceable pants graph.

\begin{definition}
We define the \textit{nondisplaceable pants graph} as the graph whose vertices correspond to isotopy classes of nondisplaceable pairs of pants, $P$, such that the complementary components of $\Sigma \ssm P$ which are not accumulated by genus have genus zero. We place an edge between two vertices if the geometric intersection of the two underlying pairs of pants equals 2.
\end{definition}

By the proof of Lemma~\ref{lem:nondisplaceableimpliesalphabeta}, a nondisplaceable pair of pants contains two separating curves, which necessitates that the third curve also 
be separating. Thus, if a curve $c_1$ on $P_1$ 
intersects a curve $c_2$ on $P_2$, then it must 
do so an even number of times. In particular, the edge condition in the definition of the nondisplaceable pants graph can be restated as follows:

\begin{corollary}\label{cor:NondisplaceableAlphaBeta}
    Let $\Sigma$ be the broom, brush, or tripod surfaces, and let $\calP(\Sigma)$ be the nondisplaceable pants graph. Then $P_1, P_2 \in \calP(\Sigma)$ are adjacent if and only if 
    there exists a unique $c_i\in P_1$ intersecting 
    exactly one curve in $P_2$ (and vice-versa).
\end{corollary}

\subsection{Adjacency properties in the nondisplaceable pants graph}\label{Subsec:AdjacencyPropertiesNDPants}

Using the interpretation of the nondisplaceable pants graph from Corollary~\ref{cor:NondisplaceableAlphaBeta} using $\alpha$ and $\beta$-curves, we wish to further understand which properties hold when two pairs of pants are adjacent in the graph.
We begin by proving that adjacent pairs of pants must share a boundary curve.

\begin{lemma}\label{lem:adjacentpantssharecurve}
Let $\Sigma$ denote either the brush, broom, or tripod surfaces and suppose that $P_1, P_2 \subset \Sigma$ are pairs of nondisplaceable pairs
pants whose boundaries consist of separating curves. If $|\del P_1 \cap \del P_2| \le 2$, then
$P_1$ and $P_2$ must share a curve.
\end{lemma}

\begin{proof}
Suppose first that $|\del P_1 \cap \del P_2| = 0$ and that $P_1$ and $P_2$ do not share a curve. Furthermore, let $\del P_i = a_i \cup b_i \cup c_i$, where $a_i$ denotes the $\alpha$-curve and $b_i$ the $\beta$-curve of $P_i$.
We will show that the $\alpha$ curves of $P_1$ and $P_2$ are shared.
In this case, without loss of generality, suppose that $a_2 \subset S_{a_1}$ where $S_{a_1}$ is the connected component of $\Sigma \ssm a_1$ consisting of none of the other curves of $P_1$. 
Then notice that there must exist a second curve of $P_2$ lying in $S_{a_1}$ in order for $P_2$ to be both disjoint from $P_1$, and a nondisplaceable pair of pants. This second curve can be constructed as the regular neighborhood of $a_1$ and $a_2$ connected by an arc in $S_{a_1}$.
Notice that this curve is neither of $\alpha$-type nor of $\beta$-type, and therefore the final curve of $P_2$ is unable to be a $\beta$-curve.
Indeed, if $\Sigma$ was either a broom or a brush, the final curve would contain all of $E_{\Gc}$, and if $\Sigma$ were a tripod the final curve would not separate $\Gb$ from $\{\Ga,\Gc\}$.
By Lemma~\ref{lem:nondisplaceableimpliesalphabeta}, this implies $P_2$ was not nondisplaceable, and thus $a_1 = a_2$.

We now let $\del P_i = a_i \cup b_i \cup c_i$ for $i \in \{1,2\}$, and without loss of generality, we assume that $a_1$ and $a_2$ intersect. By the claim above, they must intersect twice, and we must have that $a_i$ is disjoint from all $b_i$'s and $c_i$'s. 

An $\eps$-thickening of $a_1\cup a_2$ induces a subsurface $S \subset \Sigma$
which is homeomorphic to a four-times punctured sphere. We claim that
$b_1,b_2,c_1,c_2$ are each (homotopic to) curves in $\del S$. Suppose that $b_1$
is not a boundary curve of $S$, and let $R$ be the connected component of
$\Sigma \ssm S$ that contains $b_1$. Then by assumption, $b_1$ is not $\del S
\cap \del R$. Then there must be some simple closed curve, $t$, lying in $R$
disjoint from $b_1$. let $t'$ be the boundary of $S$ which lies on the same side
of $a_1$ as $R$, but which does not lie in $R$. Then $c_2$ must separate $t$ and
$t'$ from $b_1$ which forces $c_2$ to intersect $a_2$, a contradiction.
Symmetric arguments can be applied to $c_1, b_2,$ and $c_2$.

This proves that $P_1$ and $P_2$ must both lie in $S$. Since $S$ is a four-times punctured sphere induced by the $\eps$-thickening of $a_1\cup a_2$, this forces $P_1$ and $P_2$ to share a curve.
\end{proof}

Building on the previous lemma, if two pairs of nondisplaceable pants $P_1, P_2 \in \calP(\Sigma)$ are adjacent in $\calP(\Sigma)$ and the $\alpha$-curves and $\beta$-curves of two pairs of pants are pairwise disjoint, then $P_1$ and $P_2$ either share either an $\alpha$-curve or a $\beta$-curve.

\begin{lemma}\label{Lem:AdjacentPantsFacts}
Let $\Sigma$ denote either the brush, broom, or tripod surfaces and let $P_1, P_2 \in \calP(\Sigma)$ be connected by an edge. Let $\alpha_i,\beta_i\in P_i$ be the $\alpha$-type (resp. $\beta$-type) curves, and suppose that $\alpha_1 \cup \beta_1$ is disjoint from $\alpha_2 \cup \beta_2$. Then either  $\alpha_1 = \alpha_2$ or $\beta_1 = \beta_2$.
\end{lemma}

\begin{proof}
    Suppose that $\alpha_1 \cup \beta_1$ is disjoint and distinct from $\alpha_2\cup \beta_2$, and suppose that neither $\alpha_1 = \alpha_2$ nor $\beta_1 = \beta_2$.
    Without loss of generality, suppose $\alpha_1 \subset S_{\alpha_2}$. Then $\gamma_1$ must intersect $\alpha_2$.
    If $\beta_1 \subset S_{\beta_2}$, then $\gamma_1$ must intersect $\beta_2$, which implies that $P_1$ and $P_2$ are not connected by an edge in $\calP(\Sigma)$.
    If $S_{\beta_1} \supset \beta_2$, then $\gamma_2$ must intersect $\beta_1$, which again implies that $P_1$ and $P_2$ are not connected by an edge in $\calP(\Sigma)$.
\end{proof} 

\section{A uniformly hyperbolic arc and loop graph}\label{sec:arcandloopgraph}

In this subsection, we introduce a graph, denoted by $X_{\alpha,\beta}$.
The main result of this section is that the graph
$X_{\alpha,\beta}$ is uniformly $\delta$-hyperbolic.
We use this graph to define connected subgraphs of $\calP(\Sigma)$ which we use to prove that $\calP(\Sigma)$ is uniformly $\delta$-hyperbolic.

\subsection{The arc and loop graph}
In this subsection, we define the arc and loop graph, which we will denote by $X_{\alpha,\beta}$.
To begin, we introduce some notation. Fix a nondisplaceable pair of pants $P \in \calP(\Sigma)$ denote, and let $\alpha_P$ and $\beta_P$ denote the $\alpha$, and respectively $\beta$, curves in $P$.
Throughout this paper, we may slightly abuse notation and refer to pairs of pants either by their underlying surface or by their boundary curves.

Let $S_{\alpha_P,\beta_P}$ denote the subsurface of $\Sigma \ssm (\alpha_P \cup \beta_P)$ containing both $\alpha_P$ and
$\beta_P$ as boundary components, and let 
$S_{\alpha_P}$, respectively $S_{\beta_P}$, be the subsurface 
of $\Sigma \ssm (\alpha_p\cup \beta_P)$ with only $\alpha_P$, respectively $\beta_P$, as a boundary component.
Furthermore, if $\alpha$ and $\beta$ are $\alpha$ and $\beta$-curves, respectively, then let $E_{\alpha,\beta}$ denote the set of maximal ends of $S_{\alpha,\beta}$ (excluding $\alpha$ and $\beta$ as boundary components). 
If $\Sigma$ is a broom with maximal ends $\{\Ga\} \cup X$, then $E_{\alpha,\beta}$ is precisely $X \cap \overline{S_{\alpha,\beta}}$. If $\Sigma$ is a brush, then 
$E_{\alpha,\beta}$ is $E_{\Gc} \cap \overline{S_{\alpha,\beta}}$. 
If $\Sigma$ is a tripod, $E_{\alpha,\beta} = \{\Gc\}$.
We note that in the case of the broom and the brush, $E_{\alpha,\beta}$ is
homeomorphic to a Cantor set, otherwise there would exist an essential shift.

Using this notation, we now define the combinatorial graph $X_{\alpha,\beta}(\Sigma)$ where $\Sigma$ is either the broom, brush, or tripod surfaces.

\begin{definition}\label{def:alg}
    Let $\Sigma$ denote either the broom, brush, or tripod surface, and fix both an $\alpha$-curve and a $\beta$-curve in $\Sigma$, denoted by $\alpha$ and $\beta$ respectively.
    Define $X_{\alpha,\beta}(\Sigma)$ to be the combinatorial graph whose vertices consist of homotopy classes (rel $\del S_{\alpha,\beta}$) of:
    \begin{itemize}
        \item $\alpha$-loops: Loops in $S_{\alpha,\beta}$ whose endpoints lie on $\alpha$ and either separate $E_{\alpha,\beta}$, or separate $E_{\alpha,\beta}$ from $\beta$.
        \item $\beta$-loops: Loops in $S_{\alpha,\beta}$ whose endpoints lie on $\beta$ and separate $E_{\alpha,\beta}$ from $\alpha$.
        \item $\alpha - \beta$ arcs: Arcs which have one endpoint located on $\alpha$ and one endpoint on $\beta$.
    \end{itemize}
    Edges correspond to the existence of disjoint representatives between two isotopy classes of vertices. Let $\Aab \subset \Xab$ denote the subgraph spanned by all $\alpha - \beta$ arcs of $\Xab$.
\end{definition}

We note that $\Aab$ is $1$-coarsely dense in $\Xab$.
Indeed, if $x$ is an $\alpha$-loop, choose $\gamma$ to be a geodesic arc from a point $p$ on $x$ to the boundary arc $\beta$. Then following from $\alpha$ along $x$ to $p$ and then following along $\gamma$ to $\beta$ gives us an arc in $\Aab$ disjoint from $x$. The case where $x$ is a $\beta$-loop case is similar.

\subsection{Hyperbolicity of the arc and loop graph}

For the remainder of this section, we will assume that $\Sigma$ denotes either the broom, brush or tripod surfaces.
In order to prove that $X_{\alpha,\beta}(\Sigma)$ is hyperbolic, we wish to apply the Guessing Geodesics Lemma \cite{Bow2}.

\begin{lemma}[Guessing Geodesics Lemma: Proposition 3.1 of \cite{Bow} ]\label{Lem:GuessingGeodesics}
Suppose that $X$ is a connected graph, 
and that for each $x,y \in X$ there exists some 
connected subgraph $x,y \in \calU(x,y) \subset X$ 
for which there exists a uniform $M > 0$ such that:
\begin{itemize}
\item For any $x,y \in G$, if 
$d(x,y) \le 1$, then $\mathrm{diam}(\calU(x,y)) \le M$.
\item For any $x,y,z\in X$, 
$\calU(x,z) \subset \mathcal{N}_M(\calU(x,y) \cup \calU(y,z))$ (where $\mathcal{N}_M$ 
denotes the $M$-neighbourhood of a set).
\end{itemize}
Then $X$ is connected and $\delta$-hyperbolic, where 
$\delta$ depends only on $M$.
\end{lemma}

To define an appropriate connected subgraph of $\Xab(\Sigma)$, we
will introduce a modification of unicorn paths, which were introduced by
Hensel--Przytycki--Webb \cite{HPW}. Suppose that $x,y\in \Aab$ lie in minimal position, and let $e_x, e_y$ be endpoints of $x$ and $y$ lying on different curves. 
Let $x'$ (resp. $y'$) be a subarc of $x$ (resp. $y$) with endpoint $e_x$ (resp. $e_y$) with a common endpoint $\pi \in x \cap y$. If 
$x'\cup_{\pi} y'$ is an essential embedded arc, then we say that it is 
a \textit{unicorn arc} of the oriented arcs $x$ and $y$, and write $U(x,y)$ to be the set of all uncorns coming from all admissible 
orientations of $x$ and $y$.
Thus, the elements of $U(x,y)$ are precisely the arcs in $\Aab$ of the form $x' \cup y'$, where 
$x' \subset x$ and $y' \subset y$ are subarcs. We 
employ the convention that $x,y \in U(x,y)$, so 
that if $x$ and $y$ are disjoint, $U(x,y) = \{x,y\}$.

Using this definition of unicorn arcs, we show that $A_{\alpha,\beta}$ is connected.

\begin{lemma}
$\Xab$ is connected.
\end{lemma}

\begin{proof}
Remark 3.2 of \cite{HPW} shows that $U(x,y)$ are connected. Together with the observation that 
$\Aab$ is $1$-coarsely dense in $\Xab$, this immediately 
gives the result.
\end{proof}

The following lemma states that 
unicorns form thin triangles.

\begin{lemma}\label{lem:Aab_unicorns}
    Let $x,y,z\in \Aab$ be arcs. Then $U(x,y)$ is contained in the 4-neighborhood of $U(x,z)\cup U(y,z)$.
\end{lemma}

The proof of this lemma is slightly technical, and it will be proven
at the end of the section.
We note Lemma~\ref{lem:Aab_unicorns} by itself is not enough 
to show that $\Xab$ is uniformly $\delta$-hyperbolic, 
because the unicorn paths constructed are only between elements of $\Aab$. To find appropriate ``guessing geodesics'' sets in $\Xab$, we define 
$N_x = B_1(x)\cap \Aab$ for any $x\in \Xab$, and write:

\[
V(x,y) = \{z : z\in U(x',y') \text{ where } x'\in N_x, y'\in N_y\} \cup \{x,y\}
\]

Note that $V(x,y)$ is defined and is nonempty for \textit{any} elements $x,y \in \Xab$. These 
$V(x,y)$ will be the subgraphs used in the Guessing Geodesics Lemma to prove that $X_{\alpha,\beta}$ is hyperbolic.

\begin{lemma}\label{lem:GuessingGeodesicslocalcentralpants}
    Let $x,y \in \Xab$. Then $V(x,y)$ are connected subgraphs such that
    \begin{itemize}
        \item[(i)] if $d_{X_{\alpha,\beta}}(x,y) \leq 1$, then $V(x,y)$ has diameter $8$; and
        \item[(ii)] for every collection of $x,y,z \in \Xab$, $V(x,y)$ is contained in the $30$-neighborhood of $V(x,z)\cup V(y,z)$.
    \end{itemize}
\end{lemma}

\begin{proof}
We first show that the $V(x,y)$ are connected subgraphs. Fix $x_0 \in N_x$, and choose some $z \in V(x,y)$. Then $z \in U(x',y')$ for some $x'\in N_x, y'\in N_y$. In particular, 
there exists a path in $V(x,y)$ from $z$ to $y'$. But then $U(x_0,y')$ is a path from $y'$ to $x_0$, meaning that we've found a path from $z$ to $x_0$.

Next, we show that for disjoint $x,y$, $V(x,y)$ has 
diameter at most $8$. Recall that for any loop, there exists a disjoint $\alpha-\beta$ arc.
Therefore, without loss of generality, suppose that 
$x$ is an arc, and let $x'\in N_x$ and $y'\in N_y$ be arbitrary.
By Lemma~\ref{lem:Aab_unicorns}, 
$U(x',y') \subset \calN_4(U(x',x)\cup U(y',x)) 
\subset \calN_5(U(y',x))$, where the last containment 
follows because $U(x',x)$ is a single edge. But 
now $y$ is disjoint from both $y'$ and $x$, so 
$U(y',x)$ has diameter $2$, and $U(x',y') \subset B_7(x)$. Applying the definition of $V(x,y)$ and our assumption that $x$ is an arc, this proves (i).

We now prove (ii). Note that for any 
$x',x''\in N_x$, $U(x',x'')$ has diameter at most $2$ (since every surgery of $x',x''$ is disjoint from $x$). 
For $x,y\in \Xab$ and $y',y''\in N_y$, $x',x''\in N_x$, Lemma~\ref{lem:Aab_unicorns} applied to the triple $(y',x',x'')$, gives us that $U(x',y') \subset \calN_6(U(x'',y'))$. Applying the lemma 
again to the triple $(x'',y',y'')$ we see that
$U(x',y') \subset \calN_{12}(U(x'',y''))$. 

Thus, we see that $V(x,y) \subset \calN_{12}(U(x',y'))$ 
for any $x'\in N_x$ and $y' \in N_y$. In particular, 
an immediate application of Lemma~\ref{lem:Aab_unicorns} tells us that for any $x,y,z\in \Xab$, 
$V(x,z) \subset \calN_{12+12+4}(V(x,y)\cup V(y,z))$, 
as desired.
\end{proof}

Now an application of Lemma~\ref{Lem:GuessingGeodesics} allows us to conclude that $X_{\alpha,\beta}$ is uniformly hyperbolic, where the hyperbolicity constant is independent of the surface.

\begin{corollary}
    Suppose that $\Sigma$ is either the broom, brush or tripod surface.
    Then $X_{\alpha,\beta}(\Sigma)$ is uniformly $\delta$-hyperbolic, where the hyperbolicity constant is independent of the underlying surface.
\end{corollary}

All that's left is to prove Lemma~\ref{lem:Aab_unicorns}. We will need the following definition:

\begin{definition} (Access surfaces and points)
    Let $x,y \in X_{\alpha,\beta}(\Sigma)$ be two $\alpha-\beta$ arcs.
    \begin{itemize}
    \item  We say that a connected component of $S_{\alpha,\beta} \ssm (x \cup y)$ is
    an \textit{access surface} if it's intersection with $E_{\alpha,\beta}$ is nonempty.
    \item If $A$ is an access surface, 
    and $x$ is oriented, we define the \textit{access point of $A$ with respect to $y$} as the first point along $x$ that lies in $\overline{A}$. We also say that this point is an \textit{access point of $x$ with respect to $y$}. Note that 
    this access point lies in $x\cap y$. 
    \end{itemize}
\end{definition}

\begin{proof}[Proof of Lemma~\ref{lem:Aab_unicorns}]
    Let $x,y,z \in \Aab$ be $\alpha-\beta$ arcs, and 
    let $x_1\cup y_1 \in U(x,y)$ be a unicorn arc. Without loss of generality, suppose that the endpoint of $x_1$ lies on $\alpha$ and the endpoint of $y_1$ lies on $\beta$. Let $z_{\alpha}$ (resp. 
    $z_{\beta}$) be the endpoint of $z$ lying on 
    $\alpha$ (resp. $\beta$). We first orient $z$ to start at $z_{\alpha}$ and that the first point of intersection, $p$, with $x_1 \cup y_1$ lies on $y_1$. Then by following $z$ from 
    $\alpha$ to $p$ then $y_1$ to $\beta$, we get 
    an arc in $U(y,z)$ disjoint from $x_1\cup y_1$, 
    as desired. Similarly, if we orient $z$ to start at $z_{\beta}$ and that the first point of intersection, $p$, with $x_1 \cup y_1$ lies on $x_1$, 
    then we follow $z$ from $\beta$ to $p$ then $x_1$ to 
    $\alpha$, giving an arc in $U(x,z)$ disjoint from $x_1\cup y_1$.
    In either of these cases, we have that 
    $x_1\cup y_1 \in \calN_1(U(x,z)\cup U(y,z))$.

    We now suppose that when orienting $z$ to start at $z_{\alpha}$, the first point of intersection with $x_1 \cup y_1$ lies on $x_1$, and when following $z$ starting at $z_{\beta}$, the first point of intersection with $x_1 \cup y_1$ lies on $y_1$. In this case, the surgeries with $x_1\cup y_1$ as 
    above do not immediately provide an obvious unicorn 
    in $U(z,y)$ or $U(z,x)$, so we construct a loop, 
    $\delta$, that will provide a shortcut.
    We orient $x_1 \cup y_1$ to start at $\alpha$ and end at $\beta$, and follow along $x_1 \cup y_1$ until we reach an access point, $p$ of $x_1 \cup y_1$ with respect to $z$. 

    \textbf{Case 1 $p \in x_1$:} Suppose that the first access point of $x_1
    \cup y_1$ with respect to $z$ lies on $x_1$. Denote the corresponding access
    surface by $A$, and  let $X_A = E_{\alpha,\beta} \cap A$ be the set of 
    ends in $E_{\alpha,\beta}$ in $A$. Let $c_A\subset A$ be a simple closed curve separating $X_A$ from the rest of $A$, and let $\zeta\subset A$ be a simple arc from $p$ to $c_A$ whose interior 
    is disjoint from $c_A$. We define an $\alpha$-loop by starting at $\alpha$, fellow traveling $x_1$ to $p$, fellow traveling $\zeta$ to $c_A$, looping around $c_A$ back to $\zeta$, fellow traveling along $\zeta$ back to $p$, and fellow traveling along $x_1$ back to $\alpha$.
    We denote this $\alpha$ loop by $\delta$. By construction, $\delta$ can be perturbed by isotopy, if necessary, to be entirely disjoint from $x_1 \cup y_1$.
    
    Now orient $z$ to start at $\beta$ and end at $\alpha$, and let $p'$ denote the first point of intersection of $z$ with $\delta$, and let 
    $z_1\cup \delta_1$ be the unicorn surgery arising from the intersection at $p'$. By construction, 
    $z_1\cup \delta_1$ is an $\alpha-\beta$ arc 
    disjoint from $\delta$. Moreover, by the construction 
    of $\delta$, $z$ can only intersect $\delta$ along $x_1$, so actually, $\delta_1 \subset x_1$ and 
    $z_1 \cup \delta_1 \in U(z,x)$. Thus, 
    we have found the path $x_1 \cup y_1 - \delta - z_1 \cup \delta_1 \in U(x,z)$, which is a path of length 2.

    \textbf{Case 2 $p\in y_1$:} 
    We now suppose that the first access point of $x_1 \cup y_1$ with respect to $z$ lies on $y_1$. In 
    particular, this means that all access points of $x_1\cup y_1$ with respect to $z$ lie on $y_1$. 
    In this case, we will construct a $\beta$-loop, which we denote by $\delta$, that will provide the appropriate shortcut. 
    
    To construct this $\beta$ loop, begin by orienting $x_1 \cup y_1$ to start at the point $\beta$ and end at the point $\alpha$, and label each access point of 
    $x_1\cup y_1$ with respect to $z$ by $p_i$. For each access point, denote the corresponding access surface by $A_i$. For each access surface, let $c_i \subset A_i$ be a simple closed curve separating 
    all of $E_{\alpha,\beta}\cap A_i$ from $\del A_i$. Let $\zeta_i \subset A_i$ be an arc from $p_i$ to $c_i$ whose interior is disjoint from $c_i$. 
    Let $\gamma_i$ be the loop starting at $p_i$, fellow traveling along $\zeta_i$ to $c_i$, looping around $c_i$, and returning to $p_i$ along $\zeta_i$. 

    Let $y'$ be the boundary of an $\eps$-regular neighbourhood of $y_1$ (so that $y'$ forms a loop 
    that fellow travels $y_1$ to the left and right and at $p$, it loops back).
    We can find $\eps$ sufficiently small so that for each $i$, $y'\cap
    \gamma_i$ is a single point, and then we define $\delta$ to be the surgery
    of $y'$ with all of the $\gamma_i$'s. That is, we go along $y'$ until we hit
    some $\gamma_i$, then loop around $\gamma_i$ before continuing. Now, note
    that $E_{\alpha,\beta}$ is contained in the union of all access surfaces,
    meaning that the loop 
    $\delta$ separates $\alpha$ from all of $E_{\alpha,\beta}$, and hence $\delta$ is a $\beta$-loop. By construction, $\delta$ intersects 
    $x_1\cup y_1$ exactly once, near $p$, at some point $p'$. Thus, the surgery of $\delta$ and $x_1\cup y_1$ (oriented to start at $\alpha$) at $p'$ gives an $\alpha-\beta$ arc disjoint from both, and thus, 
    $d_{\Xab}(x_1\cup y_1, \delta) \le 2$. 

    We finish our proof by constructing a path between the loop $\delta$, constructed above, and an element of $U(y,z)$.
    Orient $z$ to start at $\alpha$ and end at $\beta$.
    By construction of $\delta$, the only points of intersection of $z$ with $\delta$ are on $y_1$. 
    By starting at $\alpha$ and fellow traveling $z$ until we hit the first
    point $q$ of intersection with $\delta$, $q$, we construct the
    $\alpha-\beta$ arc $z_1 \cup \delta_1$ which is disjoint from $\delta$. 

    As constructed, $\delta$ is disjoint from $y_1$, and 
    up to homotopy, $\delta \subset y_1 \cup (\cup_i \gamma_i)$, where each of the 
    $\gamma_i$'s are disjoint from $z$ (seeing as they are contained in the access surfaces). Thus, 
    making $\eps$ sufficiently small, 
    the intersection point $q$ of $\delta$ and $z$ 
    corresponds to a point of intersection, $q'$ 
    between $z$ and $y_1$. Then 
    the surgery $z_1' \cup_{q'} y_1'$ between 
    $z$ and $y_1$ is a simple arc, $\delta'$ disjoint from $z_1\cup \delta_1$. Thus, we get the path $x_1 \cup y_1 - \delta' - \delta - \delta_1\cup z_1 - z_1' \cup_{q'} y_1' 
    \in U(z,y)$, completing the proof of the claim.
\end{proof}

\section{Hyperbolicity of the nondisplaceable pants graph}

In this section, we will prove that the nondisplaceable pants graph of the broom, brush, and tripod surfaces is uniformly hyperbolic.
This proof will again make use of the Guessing Geodesics Lemma.
To apply the Guessing Geodesics Lemma, we will construct subgraphs of $\calP(\Sigma)$ by using $X_{\alpha,\beta}$.
Throughout this section, let $\Sigma$ be either a broom, brush, or tripod surface. For
$P\in \calP(\Sigma)$, let $\alpha_P$ (respectively $\beta_P$) denote the $\alpha$-curve (respectively
$\beta$-curve) of $P$. 

\subsection{Localization}
In this subsection, we describe a procedure, called localization, which allows us to construct a subset of $\calP(\Sigma)$ for any $P \in \calP(\Sigma)$.

Let $P \in \calP(\Sigma)$, and let $\alpha_P$ and $\beta_P$ denote the $\alpha$-curve, respectively $\beta$-curve, of $P$.
Let $S_{\alpha_P,\beta_P}$ denote the subsurface of $\Sigma \ssm (\alpha_P \cup \beta_P)$ containing both $\alpha_P$ and
$\beta_P$ as boundary components, and let 
$S_{\alpha_P}$ (respectively $S_{\beta_P}$) denote the subsurface 
of $\Sigma \ssm (\alpha_p\cup \beta_P)$ containing only $\alpha_P$ (respectively $\beta_P$) as a boundary component.
Let $\alpha \subset S_{\alpha_P}$ be an $\alpha$-curve, and let $\beta \subset S_{\beta_P}$ be a $\beta$-curve. In this case, we say that the pair $(\alpha,\beta)$ \textit{localizes} $P$. 

Fix a hyperbolic metric on $\Sigma$, and notice that there exists a unique geodesic arc, $x_P$, contained in $P$ connecting $\alpha_P$ and $\beta_P$. We call $x_P$ the \textit{core arc} of $P$. 
Notice that a regular neighborhood of $\alpha_P \cup \beta_P \cup x_P$ is ambiently homeomorphic to $P$.
Let $(\alpha,\beta)$ be a localizing pair of $P$ and extend $x_P$ bi-directionally obtain a simple arc with endpoints lying on $\alpha$ and $\beta$, denoted $\hat x_{\alpha\beta}$.
Notice that the arc $\hat x_{\alpha\beta}$ is unique only when the localizing pair is $(\alpha_P,\beta_P)$. 
The regular neighborhood of $\alpha \cup \hat x_{\alpha\beta}\cup \beta$ defines a pair of pants $\hat P$ which is said to be obtained from $P$ via \textit{localization by $(\alpha,\beta)$}.
Let $L_{\alpha,\beta}(P) \subset \calP(\Sigma)$ denote the set of all pants obtained from $P$ via localization by $(\alpha,\beta)$.
Define $\calL(P) = \cup_{(\alpha,\beta)} L_{\alpha\beta}(P)$, where the union is taken over all localizing pairs of $P$. 
We call $\calL(P)$ the \textit{localizing set} of $P$.
\begin{remark}
We note a few facts regarding localization:
\begin{itemize}
\item  Trivially, $P \in \calL(P)$. 
\item Unless $\alpha = \alpha_P$ and $\beta = \beta_P$, the set $L_{\alpha,\beta}(P)$ is infinite. 
\item Suppose that $P$ is a pair of pants localized by the pair $(\alpha,\beta)$. 
Since a pair of pants has no non-boundary simple closed curves, it follows that any arc $\hat x \in X_{\alpha,\beta}$ that intersects $P$ exactly twice defines a localization of $P$.
\item Let $A$ denote a finite subcollection of pants in $\calP(\Sigma)$. Then there exist curves
$\alpha$, $\beta$ that localize every $P \in A$. Indeed, such curves exist because
$\cup_{P\in A}P$ is compact, which implies that it's possible to find the appropriate curves which both separate the required ends and are disjoint from each $P \in A$.
\end{itemize}
\end{remark}

\subsection{Localization properties}
In this subsection, we will prove the properties of localization which will aid us in proving the hyperbolicity of the nondisplaceable pants graph.
The first property is that for any $P \in \calP(\Sigma)$, any pair of pants in the localizing set of $P$ is at most distance two away from $P$ in $\calP(\Sigma)$.

\begin{lemma}\label{Lem:DistanceBetweenPAndLocalizationP}
Suppose that $(\alpha,\beta)$ is a localizing pair for $P \in \calP(\Sigma)$, and let $\hat P \in
L_{\alpha\beta}(P)$. Then $d_{\calP}(P, \hat P) \le 2$.
\end{lemma}

\begin{proof}
First consider the localizing pair $(\alpha_P, \beta)$, respectively, $(\alpha,\beta_P)$, for $P$, and let $Q\in L_{\alpha_P,\beta}(P)$, respectively, $Q \in L_{\alpha,\beta_P}$.
By the localization procedure, notice that $Q$ intersects $P$ exactly twice on $\beta_P$ (resp.
$\alpha_P$), implying that $Q$ and $P$ are adjacent in $\calP(\Sigma)$. 

We may view the localization procedure to obtain $\hat P$ as a two-step localization process.
Indeed, let $\hat{x}_{\alpha,\beta}$ denote the core arc of $\hat{P}$, and let $x_P$ denote the core arc of $P$.
Begin by
localizing $P$ to $(\alpha_P, \beta_{\hat P})$ to obtain $P'$ whose core arc is the subsurface projection of $\hat{x}_{\alpha,\beta}$ to $S_{\alpha_P, \beta_{\hat P}}$. Then complete this procedure by localizing $P'$ to $(\alpha_{\hat P},\beta_{\hat P})$ to obtain $\hat{P}$. By the previous paragraph, we see that $d_{\calP}(P, \hat P) \le 2$.
\end{proof}

As a corollary of Lemma~\ref{Lem:DistanceBetweenPAndLocalizationP}, we find that $\calL(P)$ is connected and has finite diameter.

\begin{corollary}
$\calL(P)$ is connected and has diameter at most $4$.
\end{corollary}

Next, we show that if two pairs of pants, $P_1, P_2 \in \calP(\Sigma)$ are adjacent in $\calP(\Sigma)$, then if $(\alpha,\beta)$ is a localizing pair for both, the sets $L_{\alpha,\beta}(P_1)$ and $L_{\alpha,\beta}(P_2)$ have nonempty intersection.

\begin{lemma}\label{Lem:AdjacencyImpliesIntersectionLocalizations}
    Let $P_1, P_2 \in \calP(\Sigma)$, and suppose that $d_\calP(P_1, P_2) \leq 1$.
    If $(\alpha,\beta)$ is a localizing pair for both $P_1$ and $P_2$, then $L_{\alpha \beta} (P_1) \cap L_{\alpha \beta} (P_2) \neq \emptyset$.
\end{lemma}

\begin{proof}
    In this proof, we let $\alpha_i$ (resp. $\beta_i$)
    denote the $\alpha$- (resp. $\beta$-) curve of $P_i$, and let 
    $\gamma_i$ denote the curve of $P_i$ distinct from 
    $\alpha_i$ and $\beta_i$.
    Recall from Lemma~\ref{lem:adjacentpantssharecurve} that since $d_\calP(P_1, P_2) \leq 1$, $P_1$ and $P_2$ must share a curve.
    We will work through this proof based on the intersection patten of $P_1$ with $P_2$. 
    For each case, we will construct an arc, denoted $\hat x$, with endpoints lying on 
    $\alpha$ and $\beta$ that intersects $P_i$ exactly twice.

    Suppose first that $\alpha_1 = \alpha_2$ and, without loss of generality, suppose $\gamma_1$ intersects $\beta_2$.
    In this case, $\beta_1 \subset S_{\beta_2}$, $\gamma_2 \subset S_{\gamma_1}$, and there exists an arc, denoted 
    $x_1$, contained in the subsurface $P_1 \cap S_{\beta_2}$ with endpoints lying on $\beta_1$ and $\beta_2$. By construction, this arc is disjoint from $\gamma_1$. 
    Similarly, there exists an arc, denoted $x_2$, contained in the subsurface $P_2 \cap S_{\gamma_1}$ with endpoints lying on $\beta_2$ and $\alpha_1 = \alpha_2$ disjoint from $\gamma_1$ (and hence from $\gamma_2$). 
    We connect $x_1$ to $x_2$ by fellow traveling $\beta_2$ along the boundary of the subsurface $P_2 \cap S_{\beta_2}$, and denote this arc by $x$.
    By construction, $x$ intersects each of the curves $\beta_1,\beta_2, \alpha_1$ and $\alpha_2$ exactly once each, and hence, we can extend it to a simple arc $\hat x$ in $X_{\alpha,\beta}$ whose regular neighbourhood is a localization of both $P_1$ and $P_2$.

    We now suppose that $\alpha_1 = \alpha_2$, and the curves which intersect exactly twice are $\beta_1$ and $\beta_2$.
    In this case, consider the connected component of $S' = S_{\alpha_1,\beta_1} \cap S_{\alpha_2,\beta_2}$ containing $\alpha_1=\alpha_2$.
    Notice that there exists a point on the boundary of $S'$ which lies on both $\beta_1$ and $\beta_2$. Any arc $\hat{\alpha}$ connecting this point on the boundary to $\alpha_1=\alpha_2$ can be extended to a simple arc $\hat x$ in $X_{\alpha,\beta}$ whose regular neighbourhood is a localization of both $P_1$ and $P_2$. Therefore, $L_{\alpha \beta} (P_1) \cap L_{\alpha \beta} (P_2) \neq \emptyset$.

    We now suppose that $\alpha_1 = \alpha_2$, and the curves which intersect exactly twice are $\gamma_1$ and $\gamma_2$.
    Since $\beta_1$ and $\beta_2$ are disjoint, notice that without loss of generality, we may assume that $\beta_2 \subset S_{\beta_1}$.
    In this case, construct an arc in the connected component of $P_1 \cap P_2$ containing $\beta_1$ from $\beta_1$ to $\alpha_1 = \alpha_2$.
    Extend this arc, denoted $\hat{\alpha}$ to $\beta_2$ in the subsurface $S_{\beta_1} \ssm S_{\beta_2}$.
    Notice that $\hat{\alpha}$ can be extended to a simple arc $\hat x$ in $X_{\alpha,\beta}$ whose regular neighbourhood is a localization of both $P_1$ and $P_2$. Therefore, $L_{\alpha \beta} (P_1) \cap L_{\alpha \beta} (P_2) \neq \emptyset$.
    
    The case where $\beta_1 = \beta_2$ and either $\gamma_1$ intersects $\beta_2$, $\alpha_1$ intersects $\alpha_2$, or $\gamma_1$ intersects $\gamma_2$ all follow a symmetric argument to one of the above three cases.

    Now suppose that $\gamma_1 = \gamma_2$, and, without loss of generality, assume that $\alpha_2$ intersects $\beta_1$.
    This implies that 
    $\alpha_1 \subset S_{\alpha_2}$ and $\beta_2 \subset S_{\beta_1}$.
    By following a symmetric procedure to the previous paragraph, we construct an arc, denoted $y$ with endpoints lying on $\alpha_1$ and $\gamma_1 = \gamma_2$ which intersects $\alpha_2$ exactly once and is disjoint from $\beta_i$. 
    Similarly, by following a symmetric procedure to the previous paragraph, one can construct an arc, denoted $z$, with endpoints lying on $\beta_2$ and $\gamma_1 = \gamma_2$ which intersects $\beta_1$ exactly once and is disjoint from $\alpha_i$.
    We now take the unicorn arc $x = y' \cup z'$ where $y' \subset y$ is a subarc starting at $\alpha_1$ and $z' \subset z$ is a subarc ending at $\beta_2$, which is a simple arc from $\alpha_1$ to $\beta_2$.
    By construction, we have that $x$ intersects each $\alpha_2$ and $\beta_1$ at most once.
    This implies that we can extend $\hat{x}$ to a simple arc in $X_{\alpha,\beta}$ whose regular neighbourhood is a localization of both $P_1$ and $P_2$, completing the proof of our claim.

    Finally, suppose that $\gamma_1 = \gamma_2$ and $\beta_1$ intersects $\beta_2$ exactly twice. By definition, this implies that $\alpha_1$ and $\alpha_2$ are disjoint. 
    Without loss of generality, suppose that $\alpha_2 \subset S_{\alpha_1}$.
    Consider the connected component of $P_1 \cap P_2$ containing $\alpha_1$ as a boundary component, and denote this surface by $S'$.
    Notice that there exists a point on the boundary of $S'$ which lies both on $\beta_1$ and on $\beta_2$, which we denote by $p$.
    Construct an arc from $p$ to $\alpha_1$ within $S'$, and extend this arc to $\beta_2$ within the subsurface $S_{\beta_1} \ssm S_{\beta_2}$, and denote this arc by $\hat{\alpha}$.
    Notice that $\hat{\alpha}$ can be extended to a simple arc $\hat x$ in $X_{\alpha,\beta}$ whose regular neighbourhood is a localization of both $P_1$ and $P_2$. Therefore, $L_{\alpha \beta} (P_1) \cap L_{\alpha \beta} (P_2) \neq \emptyset$.

    The final case where $\gamma_1 = \gamma_2$ and $\alpha_1$ intersects $\alpha_2$ exactly twice follows a symmetric argument to the previous case.
\end{proof}

\subsection{Connectivity of the nondisplaceable pants graph}
In this subsection, we will prove the connectivity of the nondisplaceable pants graph.
Fix an $\alpha$-curve, denoted $\alpha$, and a $\beta$-curve, denoted $\beta$, and let $X_{\alpha,\beta}$ denote the corresponding arc and loop graph.
Furthermore, for any $\alpha-\beta$ arc $x \in X_{\alpha,\beta}$, let $P_x$ denote
the regular neighbourhood of $x\cup\alpha\cup\beta$. If $x \in X_{\alpha,\beta}$
is an $\alpha$-loop, respectively a $\beta$-loop, then let $P_x$ be the regular
neighbourhood of $x \cup\alpha$, respectively $x \cup\beta$. It's not hard to check that for any $x \in X_{\alpha,\beta}$, $P_x \in \calP(\Sigma)$.

\begin{lemma}\label{Lem:FourLipschitz}
The map $\phi \from X_{\alpha,\beta} \to \calP(\Sigma)$ defined by $x \mapsto P_x$ is $4$-Lipschitz. 
\end{lemma}

\begin{proof}
Let $x, y \in X_{\alpha,\beta}$.
It suffices to show that that if $x$ is an $\alpha$-loop, respectively a $\beta$-loop and $y\in X_{\alpha,\beta}$ is an $\alpha-\beta$ arc such that  $x\cap y = \emptyset$, then $P_x$ and $P_y$ either intersect twice or define the same pair of pants.
Indeed, if $x$ and $y$ are both disjoint $\alpha$ loops, then $\Sigma \ssm\{x,y\}$ must contain a subsurface which contains a subarc of $\alpha$ and a subarc of $\beta$ within the boundary components.
This must be true as by definition, one component of $\Sigma \ssm \{x\}$ must contain $\beta$, and since $x,y \in X_{\alpha,\beta}$ are disjoint, the component containing $\beta$ either contains $y$ or doesn't.
If it doesn't contain $y$, we have the required subsurface, and if it does contain $y$, then the subsurface of $\Sigma \ssm y$ containing $\beta$ is the necessary subsurface.
Thus there exists an $\alpha-\beta$ arc disjoint from $x$ and $y$.
We note that the above argument holds if $x$ and $y$ are both disjoint $\beta$ loops.
Suppose now that $x,y \in X_{\alpha,\beta}$ are disjoint $\alpha-\beta$ arcs.
Since $x$ and $y$ are both compact, there exists a connected component of $S_{\alpha,\beta} \ssm \{x,y\}$
whose closure intersects $\beta$ and contains a nonempty clopen subset of $E_{\alpha,\beta}$.
Therefore, we can construct the following disjoint $\beta$-loop from both $x$ and $y$: consider a neighborhood $U$ of $E_{\alpha,\beta}$ disjoint from $x$ and $y$ and let the boundary of this neighborhood be denoted by $u$. Connect $u$ to $\beta$ via a geodesic arc, denoted by $a$, in $\Sigma \ssm\{x,y\}$. The regular neighborhood of $u \cup a$ defines the desired $\beta$-loop.
We remark that by definition, it is never the case that an $\alpha$-loop is disjoint from a $\beta$ loop.

We will begin by showing that if $x$ is an $\alpha$-loop and $y\in X_{\alpha,\beta}$ is an $\alpha-\beta$ arc such that  $x\cap y = \emptyset$, then $P_x$ and $P_y$ either intersect twice or define the same pair of pants.
We note that if $x$ were a $\beta$-loop, then the proof follows a nearly identical argument.

Therefore, suppose that $x$ is an $\alpha$-loop.
Since $y$ is an $\alpha-\beta$ arc, notice that the regular neighborhood of $\alpha \cup \beta \cup y$ will define a pair of pants $P_y = \{\alpha,\beta,\hat y\}$.
Now notice that the component of $S_{\alpha,\beta}\ssm x$ containing $y$ has two boundary components: one equivalent to $\alpha$, respectively $\beta$, and one of the form $x \cup a'$ where $a$ denotes a subarc of $\alpha$.
We denote this subsurface by $S^y$.
If this subsurface is equivalent to an annulus containing $y$, then $x$ and $y$ define the same pair of pants as the regular neighborhood of $\alpha \cup x$ containing $\beta$ is homotopic to $\beta$. 
Indeed, the regular neighborhood of $y \cup \beta$ is isotopic to the boundary component of the annulus, implying that the pants defined by $x$ and $y$ are isotopic.

On the other hand, suppose that the subsurface $S^y$ is not an annulus.
Notice that since $x$ is an $\alpha$-loop, that $\alpha$ must be a boundary curve of $P_x$, and subsequently a boundary curve of $P_y$.
Additionally, $\alpha$ and $\beta$ are both disjoint from all boundary curves of $P_x$ by construction.
Let $b$ denote the boundary curve of $S_{\alpha,\beta}\ssm x$ of the form $x \cup a'$ where $a'$ is a subarc of $\alpha$ such that $S_b \supset \beta$, and let $c$ denote the boundary curve of $S_{\alpha,\beta}\ssm x$ consisting of $x$ together with a subarc of $\alpha$ such that $S_c \subset \hat{y}$.
It remains to show that $b$ must intersect $\hat{y}$ exactly twice.
Indeed, by definition, $b$ is isotopic to the boundary of the subsurface $S^y$ which is not equivalent to ${\beta}$.
On the other hand, the regular neighborhood of $\alpha\cup\beta \cup y$ is by definition an $\epsilon$ neighborhood around $\alpha\cup\beta \cup y$. This implies the the only portion of this regular neighborhood which leaves the subsurface $S^y$ is the portion which remains within an $\epsilon$-neighborhood of $\alpha$. This implies that the regular neighborhood of $\alpha\cup\beta \cup y$ exits and reenters the subsurface $S^y$, implying that, by construction, $P_x$ and $P_y$ intersect exactly twice and are thus adjacent in $\calP(\Sigma)$.

This, together with the first paragraph, shows that the map $\phi \from X_{\alpha,\beta} \to \calP(\Sigma)$ is 4-Lipschitz.
\end{proof}

Using this 4-Lipschitz map together with the connectivity of the arc and loop graph, we prove that the nondisplaceable pants graph is connected.

\begin{lemma}\label{Lem:PConnected}
Let $\Sigma$ be a broom, brush, or tripod. Then $\calP(\Sigma)$ is connected.
\end{lemma}

\begin{proof}
Let $P_1, P_2 \in \calP(\Sigma)$ and let $(\alpha,\beta)$ be a localizing pair for both $P_1$ and $P_2$.
Let $\hat P_i \in L_{\alpha,\beta}(P_i)$, and let $\hat x_i$ be the unique arc contained
in $\hat P_i$ with endpoints lying on $\alpha_{\hat P_i} = \alpha$ and $\beta_{\hat P_i} = \beta$. 
By Lemma~\ref{lem:GuessingGeodesicslocalcentralpants}
$U(\hat x_1, \hat x_2)$ is connected and contains $\hat x_1$ and $\hat x_2$, so there exists
a path $\hat x_1 = y_0, y_1, \hdots, y_n = \hat x_2$ in
$X_{\alpha,\beta}$. 
By Lemma~\ref{Lem:FourLipschitz}, $d_{\calP}(P_{y_i},
P_{y_{i+1}}) \le 4$, implying that there is a path from $\hat P_1
= P_{y_0}$ and $\hat P_2 = P_{y_n}$.
Observing that $\hat P_i$ is distance at most $2$ from $P_i$ in $\calP(\Sigma)$ completes the proof.
\end{proof}

\subsection{Hyperbolicity of the nondisplaceable pants graph}
In this subsection, we will prove the hyperbolicity of the nondisplaceable pants graph $\calP(\Sigma)$.
To prove this result, we will again apply the Guessing Geodesics Lemma \cite{Bow}.
In this application of the Guessing Geodesics Lemma, the connected subgraphs we will use are a union of unicorn paths in in the arc and loop graph defined in Section~\ref{sec:arcandloopgraph}. We begin by defining the appropriate connected subgraph.

Let $P \in \calP(\Sigma)$, and fix a localizing pair for $P$, denoted $(\alpha,\beta)$.
By applying the map $\psi \from L_{\alpha,\beta}(P) \to X_{\alpha,\beta}$ defined by sending $P' \in L_{\alpha,\beta}(P)$ to its core arc $x_p \in X_{\alpha,\beta}$, we are able to view $L_{\alpha,\beta}(P)$ as a subset of $X_{\alpha,\beta}$.
The restriction of the map $\phi$ from Lemma~\ref{Lem:FourLipschitz} to the set of arcs $X_{\alpha,\beta}(P)$ provides us with a well-defined inverse for $\psi$ (under homotopy of endpoints on the boundary components $\alpha$ and $\beta$).
We define $$X_{\alpha,\beta}(P) = \{\psi(P') : P' \in L_{\alpha,\beta}(P)\}  \subset X_{\alpha,\beta}.$$

Recall that for $x,y \in X_{\alpha,\beta}$, we denote
the set of all unicorn paths between $x$ and $y$ by $U(x,y)$. 
Using the unicorn paths between two elements in $X_{\alpha,\beta}(\Sigma)$, we will define a connected subgraph associated to any two $P_1,P_2 \in \calP(\Sigma)$.

\begin{definition}
For $P_1, P_2 \in \calP(\Sigma)$, define
$$
\calU(P_1, P_2) = \calL(P_1) \cup \calL(P_2) \cup \bigcup_{\alpha, \beta} 
\left(
\bigcup_{x_1 \in X_{\alpha,\beta}(P_1), x_2 \in X_{\alpha,\beta}(P_2)} \{P_x | x \in U(x_1, x_2)\}
\right),
$$
where the final union is taken over all localizing pairs $(\alpha,\beta)$ which localize both $P_1$ and $P_2$ simultaneously.
\end{definition}

We now prove that the assumptions of the Guessing Geodesics lemma hold for $\calP(\Sigma)$.

\begin{lemma}\label{lem:guessinggeodesicsP}
    Let $P_1, P_2 \in \calP(\Sigma)$. Then $\calU(P_1,P_2)$ is a connected subgraph, and 
    there exists some $M$ such that 
    \begin{itemize}
        \item if $d_{\calP}(P_1,P_2) \leq 1$, then $\calU(P_1,P_2)$ has diameter $\le 19$; and
        \item for every collection of three unique $P_1, P_2, P_3 \in \calP(\Sigma)$, $\calU(P_1,P_2)$ is contained in the $18$-neighborhood of $\calU(P_1,P_3)\cup \calU(P_2,P_3)$.
    \end{itemize}
\end{lemma}

\begin{proof}
We begin by showing that $\calU(P_1,P_2)$ is a connected subgraph. Let $P, P' \in \calU(P_1,P_2)$, where $P_1,P_2 \in \calP(\Sigma)$.
Notice that $P, P'$ both lie in $\calP(\Sigma)$, and by choosing a pair $(\alpha,\beta)$ which localize both $P$ and $P'$, $(\alpha,\beta)$ also localizes $P_1$ and $P_2$. 
We may now follow the argument provided in the proof of Lemma~\ref{Lem:PConnected} to conclude that $\calU(P_1,P_2)$ is a connected subgraph.

We will now show that $\calU(P_1,P_2)$ has finite diameter.
Let $P_1, P_2 \in \calP(\Sigma)$ and suppose that $d_{\calP}(P_1,P_2) \leq 1$.
Suppose further that $P \in \calU(P_1,P_2)$.
We will show that $P$ is within bounded distance of either $P_1$ and $P_2$.

Since $P \in \calU(P_1,P_2)$, by definition it is the regular neighbourhood of $\alpha \cup \beta \cup x_P$ for some core curve $x_P \in U(x_1,x_2)$, where $x_1 \in X_{\alpha,\beta}(P_1)$ and $x_2 \in X_{\alpha,\beta}(P_2)$.
Let $P_3 \in L_{\alpha,\beta}(P_1)\cap L_{\alpha,\beta}(P_2)$, which exists by Lemma~\ref{Lem:AdjacencyImpliesIntersectionLocalizations}. 
By definition, $x_3 = x_{P_3}$ lies in both $X_{\alpha \beta}(P_1)$ and $X_{\alpha \beta}(P_2)$.
By Lemma~\ref{lem:Aab_unicorns}, $x_P$ is distance 4 away from some $z \in U(x_1,x_3) \cup U(x_2,x_3)$. Applying Lemma~\ref{Lem:FourLipschitz} implies that that $P_{x_P} = P$ is distance at most $16$ away from $P_z$. We claim that $P_z$ is close to either $P_1$ or $P_2$, and subsequently both $P_1$ and $P_2$.

Suppose first that $z \in U(x_1,x_3)$, which implies that $z =x_1' \sqcup x_3'$ where $x_1' \subset x_1$ and $x_3' \subset x_3$.
We claim that the regular neighborhood of $z \cup \alpha \cup \beta$ is a localization of $P_1$.
Indeed, suppose that $\alpha_{P_1}, \beta_{P_1},$ and $\gamma_{P_1}$ be the boundary components of $P_1$.
If $z$ intersects $\alpha_{P_1}$, respectively $\beta_{P_1}$, more than once then it must intersect this boundary component at least $3$ times as $\alpha \subset S_{\alpha_{P_1}}$, respectively $\beta \subset S_{\beta_{P_1}}$.
This implies that either $x_1' \subset x_1$ or $x_3' \subset x_3$ must intersect a boundary component of $\alpha_{P_1}$, respectively $\beta_{P_1}$ more than once implying that either $x_1$ or $x_3$ intersects a boundary component of $P_1$ more than once.
If $x_1$ intersects a boundary component of $P_1$ more than once, this contradicts that $x_1 \in X_{\alpha,\beta}(P_1)$.
If $x_3$ intersects a boundary component of $P_1$ more than once, this contradicts that $P_3 \in L_{\alpha,\beta}(P_1)\cap L_{\alpha,\beta}(P_2)$.
Therefore the regular neighborhood of $z \cup \alpha \cup \beta$ is a localization of $P_1$.
A symmetric argument shows that if $z \in U(x_2,x_3)$, then $z \cup \alpha \cup \beta$ is the core curve of a localization of $P_2$. 
Since $P_z$ is a localization of either $P_1$ or $P_2$, we apply Lemma~\ref{Lem:DistanceBetweenPAndLocalizationP} which implies that $P_z$ is distance two away from either $P_1$ or $P_2$.
Therefore, we have shown that any $P \in \calU(P_1,P_2)$ is distance at most 16 away from $P_z$, which is at most distance two away from either $P_1$ or $P_2$, and is this at most distance 19 away from both $P_1$ and $P_2$.

We will now complete the proof by showing that for any 
$P_1, P_2, P_3\in \calP(\Sigma)$ and
$Q \in \calU(P_1, P_2)$ that $Q \in \calN_{18}(\calU(P_1,P_3)\cup \calU(P_2,P_3))$. 

By definition, $Q$ is a pair of pants such that its core arc $x_Q$ lies in $U(x_1, x_2)$ for some $x_1, x_2 \in X_{\alpha,\beta}$ where $(\alpha,\beta)$ localizes $P_1$ and $P_2$.
If $(\alpha,\beta)$ localizes $P_3$, then for any $x_3 \in X_{\alpha,\beta}(P_3)$ there exists some $y \in U( x_1, x_3) \cup U(x_2, x_3)$ at most distance 4 away from $x_Q$ by Lemma~\ref{lem:Aab_unicorns}. 
This implies that $P(y) \in \calU(P_1,P_3)\cup \calU(P_2,P_3)$ and $P(y)$ is at most distance $16$ away from $P(x_Q) = Q$ by Lemma~\ref{Lem:FourLipschitz}.

We now suppose that $(\alpha,\beta)$ is not a localizing pair for $P_3$.
In this case, 
let $(\alpha', \beta')$ be a localizing pair for $P_1,P_2,$ and $P_3$, implying that $\alpha' \subset S_{\alpha}$ and $\beta' \subset S_{\beta}$. 
Recall that $x_Q$ is an element of the unicorn path between $x_1$ and $ x_2$ in $X_{\alpha,\beta}$.
We may extend the arcs $x_1$ and $x_2$ to $\alpha'$ and $\beta'$ in a way that preserves the intersection structure between $x_1$ and $x_2$.
Call these newly extended arcs $x_1'$ and $x_2'$ and let $U(x_1', x_2') \subset X_{\alpha',\beta'}$ be the unicorn path defined between $x_1'$ and $x_2'$.
Notice that if $x_i \in U(x_1, x_2)$, then the arc $x_i' \in U(x_1',x_2')$, which by construction is the extension of $x_i$ to $x_i'$, is at most distance 2 away by following the same argument found in Lemma~\ref{Lem:DistanceBetweenPAndLocalizationP}.
Since $\alpha'$ and $\beta'$ localize $P_1$, $P_2$ and $P_3$, then for any
$x_3 ' \in X_{\alpha',\beta'}(P_3)$ there exists some $y \in U(x_1', x_3') \cup U(x_2', x_3')$ at most distance 4 away from $x_Q'$ by Lemma~\ref{lem:Aab_unicorns}.
This implies that $P(y) \in \calU(P_1,P_3)\cup \calU(P_2,P_3)$ and $P(y)$ is at most distance $18$ away from $P(x_Q) = Q$.
\end{proof}

We may now apply Lemma~\ref{Lem:GuessingGeodesics} to conclude that $\calP(\Sigma)$ is uniformly hyperbolic, where the hyperbolicity constant is independent of the surface.

\begin{corollary}\label{Cor:HypP}
    Suppose that $\Sigma$ is either the broom, brush or tripod surface.
    Then $\calP(\Sigma)$ is uniformly $\delta$-hyperbolic, where the hyperbolicity constant is independent of the underlying surface.
\end{corollary}

\section{A quasi-isometry between the pants graph and mapping class group}
In this subsection, we will prove that there exists a quasi-isometry between the mapping class group and the nondisplaceable pants graph for the broom, brush and tripod surfaces.
This proof relies on the concept of coarsely bounded generating sets, so we will begin by giving a brief overview of the topic, before moving on to the proof.

\subsection{Coarsely Bounded Generating Sets}

We begin by introducing the concept of coarse boundedness, a useful tool when studying the geometry of the mapping class group. For more details on the connections between coarse boundedness and the geometry of the mapping class group, we recommend the work of Mann-Rafi \cite{MR}, as well as the book by Rosendal \cite{Rose}.

We remind the reader that a subset $X$ of a group $G$ is \textit{coarsely
bounded}, abbreviated CB, if every compatible left-invariant metric on $G$ gives $X$
finite diameter.
Examples of coarsely bounded sets include finite sets and compact sets. There are also two other viewpoints that are sometimes easier to work with.

\begin{definition}
    A \emph{continuous left-invariant pseudo-metric} on a topological group $G$ is a continuous map $d \colon G \times G \to \RR^{\ge 0}$ satisfying the following conditions:
    \begin{itemize}
        \item For each $g \in G$, $d(g, g) = 0$.
        \item For each $g, h \in G$, $d(g, h) = d(h, g)$.
        \item For each $g, h, k \in G$, $d(g, k) \le d(g, h) + d(h, k)$.
        \item For each $g, h, k \in G$, $d(gh, gk) = d(h, k)$.
    \end{itemize}
\end{definition}
Note that a pseudo-metric differs from a metric in that $d(g, h)$ may be zero even when $g \ne h$.
\begin{definition}
    A \emph{continuous length function} on a topological group $G$ is a continuous map $\ell \colon G \to \RR^{\ge 0}$ satisfying the following conditions:
    \begin{itemize}
        \item $\ell(e) = 0$, where $e$ is the identity element of $G$.
        \item For each $g \in G$, $\ell(g) = \ell(g^{-1})$.
        \item For each $g, h \in G$, $\ell(gh) \le \ell(g) + \ell(h)$.
    \end{itemize}
\end{definition}

We can combine these concepts via the following well-known lemma. 
\begin{lemma}
    Given a subset $A$ of a topological group $G$, the following are equivalent:
    \begin{itemize}
        \item The set $A$ is coarsely bounded in $G$.
        \item For every continuous left-invariant pseudo-metric on $G$, $A$ has bounded diameter.
        \item For every continuous length function on $G$, the elements of $A$ have bounded length.
    \end{itemize}
\end{lemma}
\begin{proof}
    There is a straightforward dictionary between continuous actions on metric spaces, continuous left-invariant pseudo-metrics, and continuous length functions. Given a pseudo-metric $d$ on $G$, $G$ acts on the metric space $G / \sim$, where $g \sim h$ whenever $d(g, h) = 0$. Likewise, if $G$ acts on a metric space $X$, we can pick a basepoint $x \in X$ and let $d(g, h) = d_X(g \cdot x, h \cdot x)$.

    For the second equivalence, consider a pseudo-metric $d$ on $G$. We can define a length function $\ell(g) = d(g, e)$ where $e$ is the identity element of $G$. Likewise, given a length function $\ell$ on $G$, we can define a pseudo-metric $d(g, h) = \ell(g^{-1}h)$.
\end{proof}

The value of coarsely bounded sets, for our purposes, lies in the following result.
\begin{proposition}[Proposition 2.52 of \cite{Rose}]\label{maximal}
    If a topological group $G$ is generated by a coarsely bounded set $A$, then the word metric $d_A$ with respect to $A$ is not only continuous and left-invariant, but maximal in the following sense: given any continuous left-invariant pseudo-metric $d$ on $G$, there exist constants $K \ge 1$ and $C \ge 0$ such that, for all $g, h \in G$,
    \[d(g, h) \le K d_A(g, h) + C.\]
\end{proposition}
Note that in particular, this implies that the quasi-isometry type of $G$ is independent of the choice of coarsely bounded generating set. A topological group $G$ admitting a coarsely bounded generating set is called \emph{CB-generated}.

For those familiar with geometric group theory, coarsely bounded generating sets play a very similar role to finite generating sets. For example, it can be shown (Lemma 5.51 of \cite{Rose}) that if $G$ is CB-generated and $H$ is a finite-index subgroup of $G$, then $H$ is quasi-isometric to $G$. Other comparisons are less fruitful. For instance, if $G$ is finitely generated, with a subgroup $H$ isomorphic to $\ZZ^2$, $G$ cannot be hyperbolic. The same does not apply when $G$ is CB-generated, since $H$ may be coarsely bounded in $G$. However, a more restrictive result is still possible.

\begin{lemma}\label{flats}
    Let $G$ be a CB-generated group, and suppose there are commuting elements $g_1, g_2 \in G$ and length functions $\ell_1, \ell_2$ on $G$ such that $\ell_1(g_1^n g_2^m)$ is proportional to $n$ and $\ell_2(g_1^n g_2^m)$ is proportional to $m$. Then $G$ is not hyperbolic with respect to any coarsely-bounded generating set.
\end{lemma}
\begin{proof}
    Let $\ell$ be a new length function defined as the sum of $\ell_1$ and $\ell_2$. Then $\{(g_1^i, g_2^i, g_1^ig_2^i)\}_{i \in \ZZ^+}$ form the vertices of triangles that are arbitrarily thick with respect to $\ell$. By Proposition~\ref{maximal}, these triangles are also arbitrarily thick with respect to the word metric on $G$ given by a coarsely bounded generating set. 
\end{proof}

\subsection{The quasi-isometry}

In this subsection, we will prove that the mapping class group and the nondisplaceable pants graph of the broom, brush and tripod surfaces are quasi-isometric.
The authors would like to thank Anschel Schaffer-Cohen for providing the proof of Proposition~\ref{prop:qi} for the broom surface, and providing the authors with the tools to prove the result for the brush and tripod surfaces.

\begin{proposition}\label{prop:qi}
    Let $\Sigma$ be a broom, brush, or tripod  surface with tame end space and CB-generated mapping class group. Then the mapping class group $\Map(\Sigma)$ is quasi-isometric to $\calP(\Sigma)$.
\end{proposition}

\begin{proof}
    Observe that $\calP(\Sigma)$ is a connected simplicial graph by Lemma~\ref{Lem:PConnected}, and that the action of $\Map(\Sigma)$ is vertex-transitive by Lemma~\ref{lem:changeofpants}. 
    It follows from Lemma 2.5 of \cite{SC} that the orbit map of this action is a quasi-isometry if and only if for some fixed vertex $P \in \calP(\Sigma)$, the set $\{g \in \Map(\Sigma) \mid d(P, g P) \le 1\}$ is coarsely bounded in $\Map(\Sigma)$.

    Fix  some $P \in \calP(\Sigma)$ as our base vertex.
    In this proof, we will begin by analyzing the stabilizer corresponding to $P$.
    Afterwards, we will analyze the additional mapping classes required to send $P$ to a vertex $P' \in \calP(\Sigma)$ such that $d_{\calP}(P,P') \leq 1$.
    For both the analysis of the stabilizer and the vertices which are at most distance 1 from $P$, we will need to do some casework and handle the broom, brush, and tripod separately.

    Let $P \in \calP(\Sigma)$ be our base vertex, and let $\calV_P = \{g \in
    \Map(\Sigma) : g|_{P} = id\}$. The stabilizer of the vertex corresponding to
    $P$ in $\calP(\Sigma)$ is simply the set of mapping classes with a representative
    that restricts to a homeomorphism of $P$. Since the pure mapping class group
    of a pair of pants is trivial, this stabilizer is almost identical to
    $\calV_P$.
    Note that $\calV_P$ is itself coarsely bounded by Theorem 1.4 of \cite{MR} since $P$ partitions the space of ends into stable, complementary components.
    We now show that the stabilizer $\Stab(P)$ is coarsely bounded in the case of the broom, brush, and tripod. A mapping class $\phi \in \Stab(P) \setminus \calV_P$ must nontrivially permute the cuffs of $P$.

		In the case of the broom, $\phi$ cannot change the cuff whose complementary
		component contains the singleton maximal end, so it must swap the cuffs
		whose complementary components partition the Cantor set of ends. Thus if
		$\rho$ is any mapping class swapping these two cuffs, $\vphi \in \calV_P
		\cup \rho \calV_P$. Since the translate of a coarsely bounded set is
		coarsely bounded, and the union of two coarsely bounded sets is coarsely
		bounded, the stabilizer $\Stab(P)$ must be coarsely bounded (and in fact is
		equal to $\calV_P \cup \rho \calV_P$).
    
		In the case of the brush, let $\calM(\Sigma) = \{\Ga,\Gb\}$, and let
		$E_\Gc$ denote the equivalence class of ends that is an immediate
		predecessor of $\Ga$ but not $\Gb$.  Any $\phi \in \Stab(P)$ cannot change
		the cuff whose complementary component contains $\Ga$ or the component
		containing $\Gb$, or the component containing some of $E_{\Gc}$. Thus,
		$\Stab(P) = \calV_P$ which is coarsely bounded.

		In the case of the tripod, let $\calM(\Sigma) = \{\Ga,\Gb,\Gc\}$. A mapping
		class $\phi \in \Stab(P)$ can permute the cuffs,
		$\alpha_P=c_1,\beta_P=c_2,\gamma_P=c_3$ of $P$ freely. Let $h_1,h_2,h_3$ be
		mapping classes such that $h_i$ swaps $c_i$ and $c_{i+1 \pmod 3}$, and let
		$r$ be an order 3 mapping class that sends $c_1$ to $c_2$, $c_2$ to $c_3$,
		and $c_3$ to $c_1$. Let $H = \{e,h_1,h_2,h_3,r,r^2\}$. Then $\phi \in
		H\calV_P$, which is coarsely-bounded as a union of coarsely-bounded sets.
        
        \begin{figure}
        \centering
        \begin{overpic}[width=0.5\textwidth]{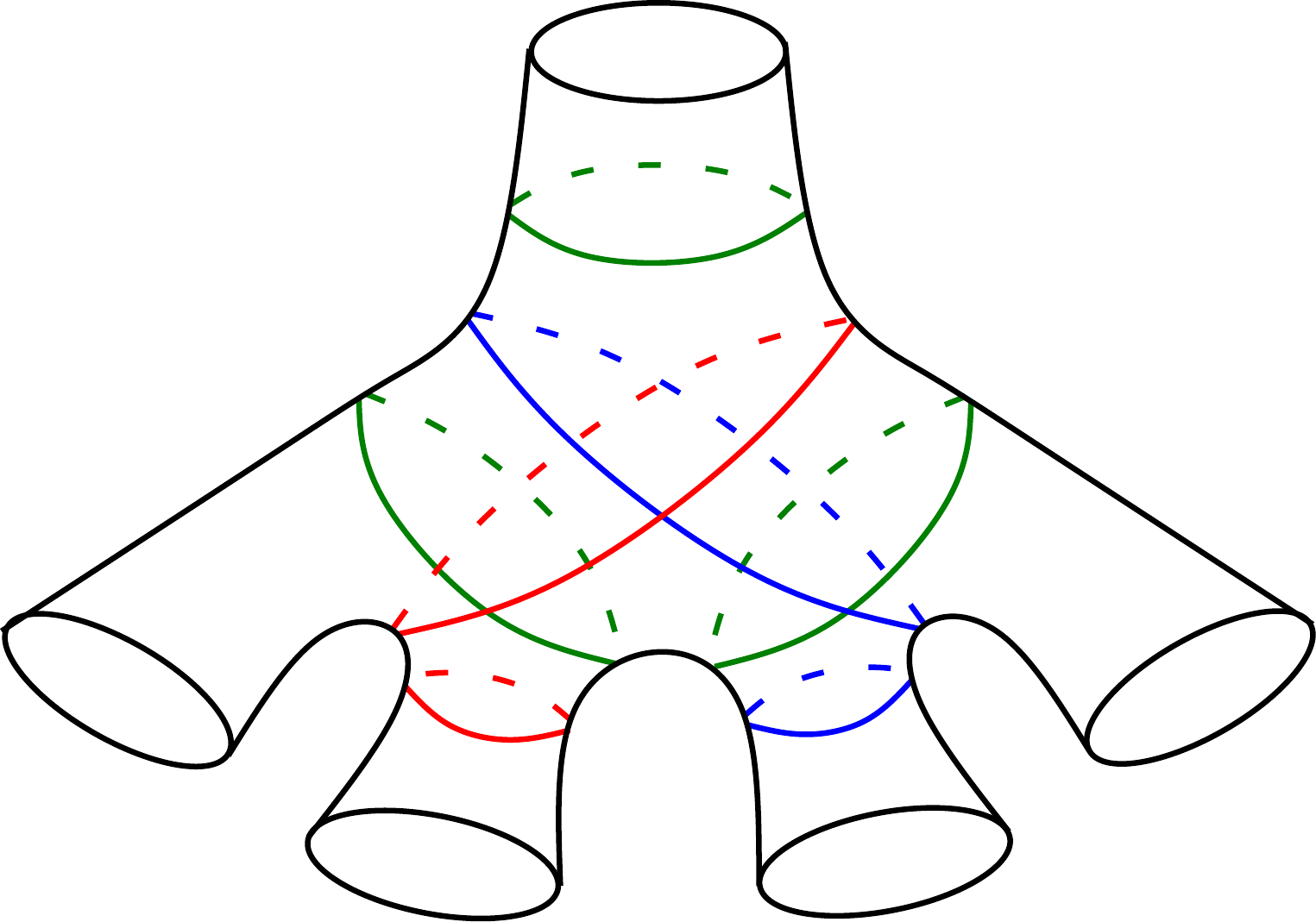}
            \put(65,55){\textcolor{ForestGreen}{$\alpha$}}
            \put(75,40){\textcolor{ForestGreen}{$\beta$}}
            \put(22,40){\textcolor{ForestGreen}{$\gamma$}}
            \put(49,72){$\Ga$}
            \put(25,-5){$C$}
            \put(70,-5){$C$}
            \put(0,10){$\calU$}
            \put(95,10){$\calU$}
        \end{overpic}
        \caption{The curves used to define the mapping classes $f_1$ and $f_2$ for both the broom and brush surfaces.}
        \label{figure:fmaps}
        \end{figure}
        
        \begin{figure}
        \centering
        \begin{overpic}[width=0.4\textwidth]{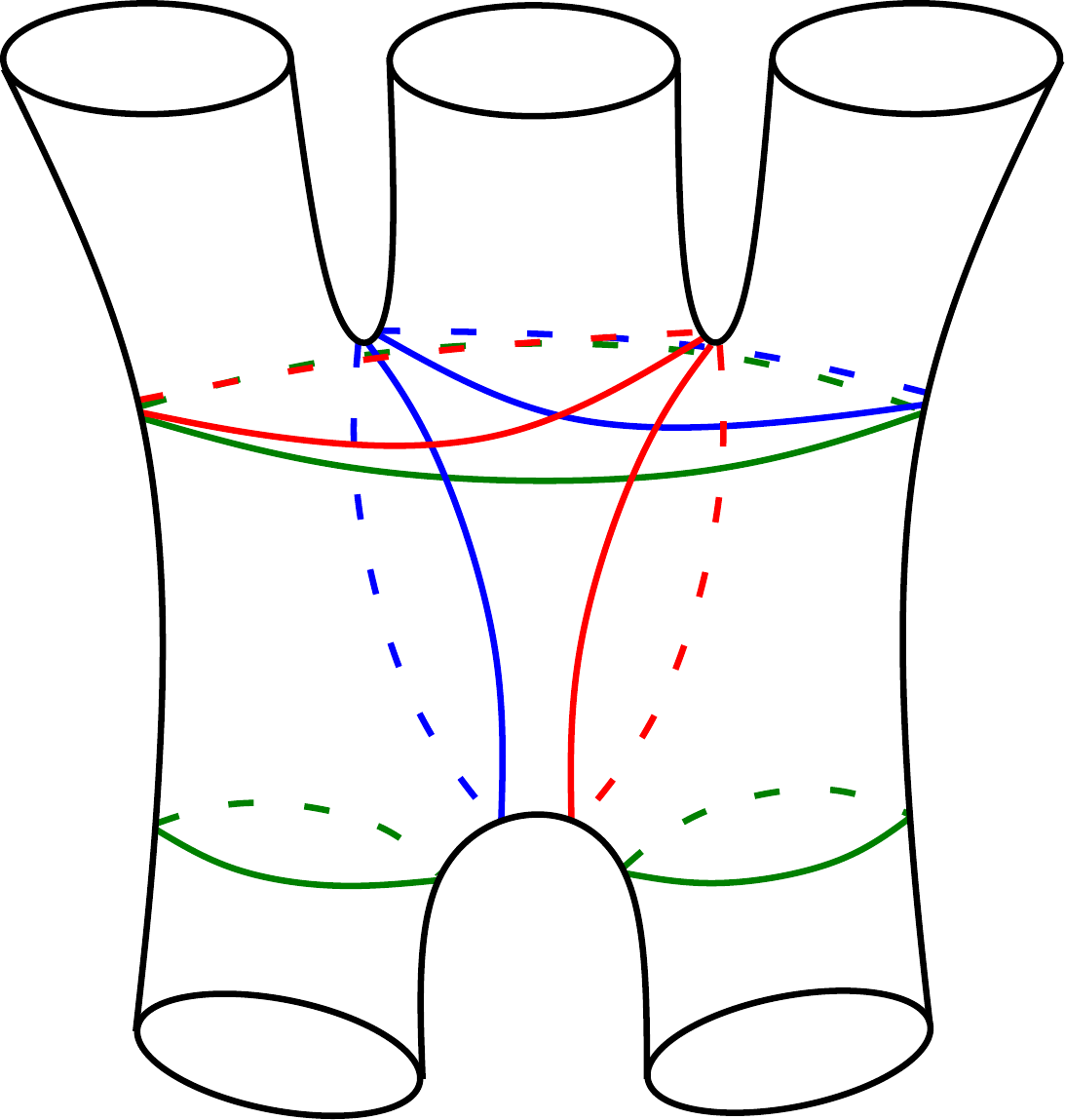}
            \put(85,62){\textcolor{ForestGreen}{$\alpha$}}
            \put(83,22){\textcolor{ForestGreen}{$\beta$}}
            \put(7,22){\textcolor{ForestGreen}{$\gamma$}}
            \put(47,102){$\Ga$}
            \put(22,-5){$C$}
            \put(70,-5){$C$}
            \put(0,99){$\calU$}
            \put(91,99){$\calU$}
        \end{overpic}
        \caption{The curves used to define the mapping classes $g_1$ and $g_2$ for both the broom and brush surfaces.}
        \label{figure:gmaps}
        \end{figure}
        
        \begin{figure}
        \centering
        \begin{overpic}[width=0.5\textwidth]{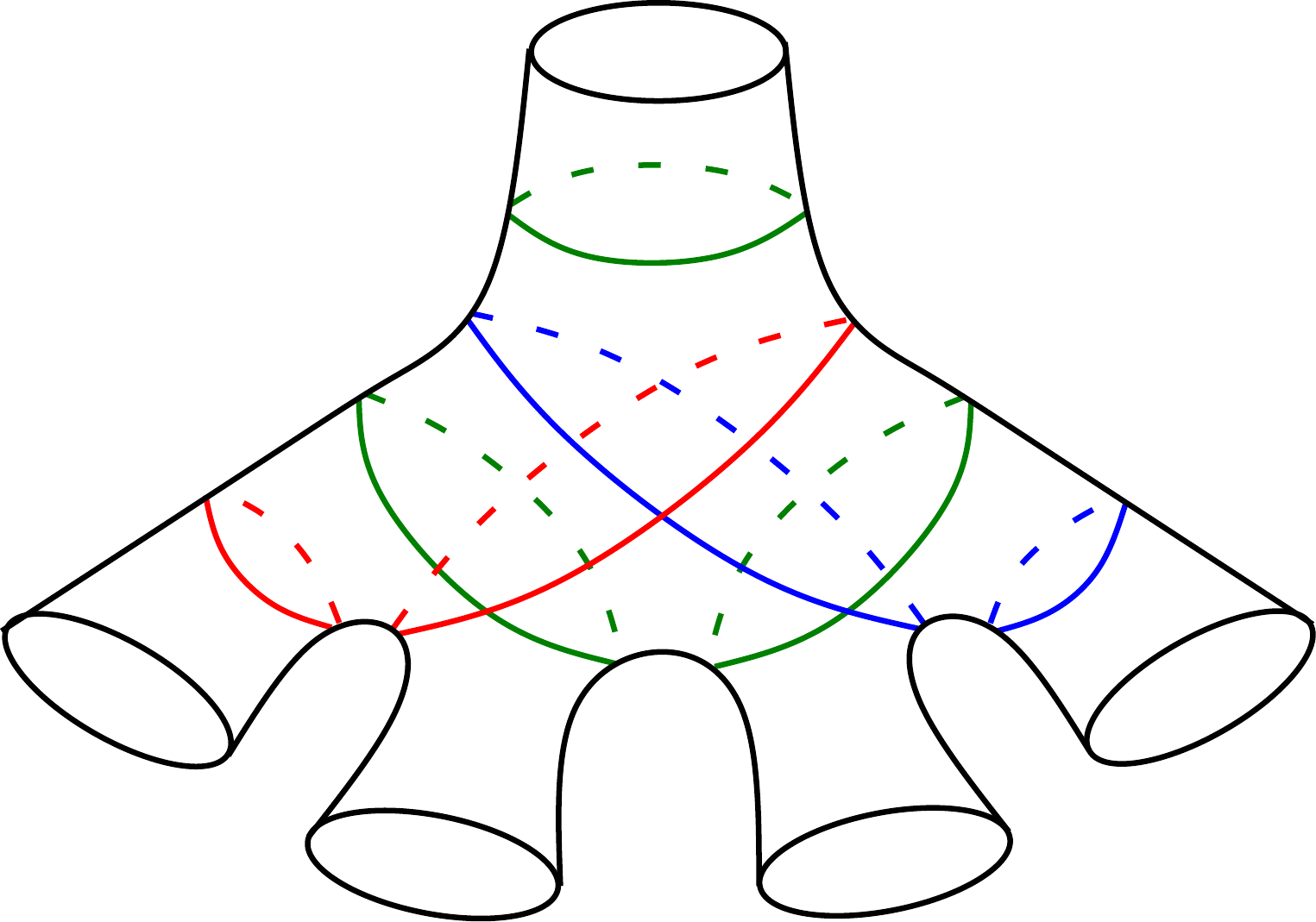}
            \put(65,55){\textcolor{ForestGreen}{$\alpha$}}
            \put(75,40){\textcolor{ForestGreen}{$\beta$}}
            \put(22,40){\textcolor{ForestGreen}{$\gamma$}}
            \put(49,72){$\Ga$}
            \put(25,-5){$\calU$}
            \put(70,-5){$\calU$}
            \put(0,10){$C$}
            \put(95,10){$C$}
        \end{overpic}
        \caption{The curves used to define the mapping classes $h_1$ and $h_2$ for both the broom surface.}
        \label{figure:hmaps}
        \end{figure}
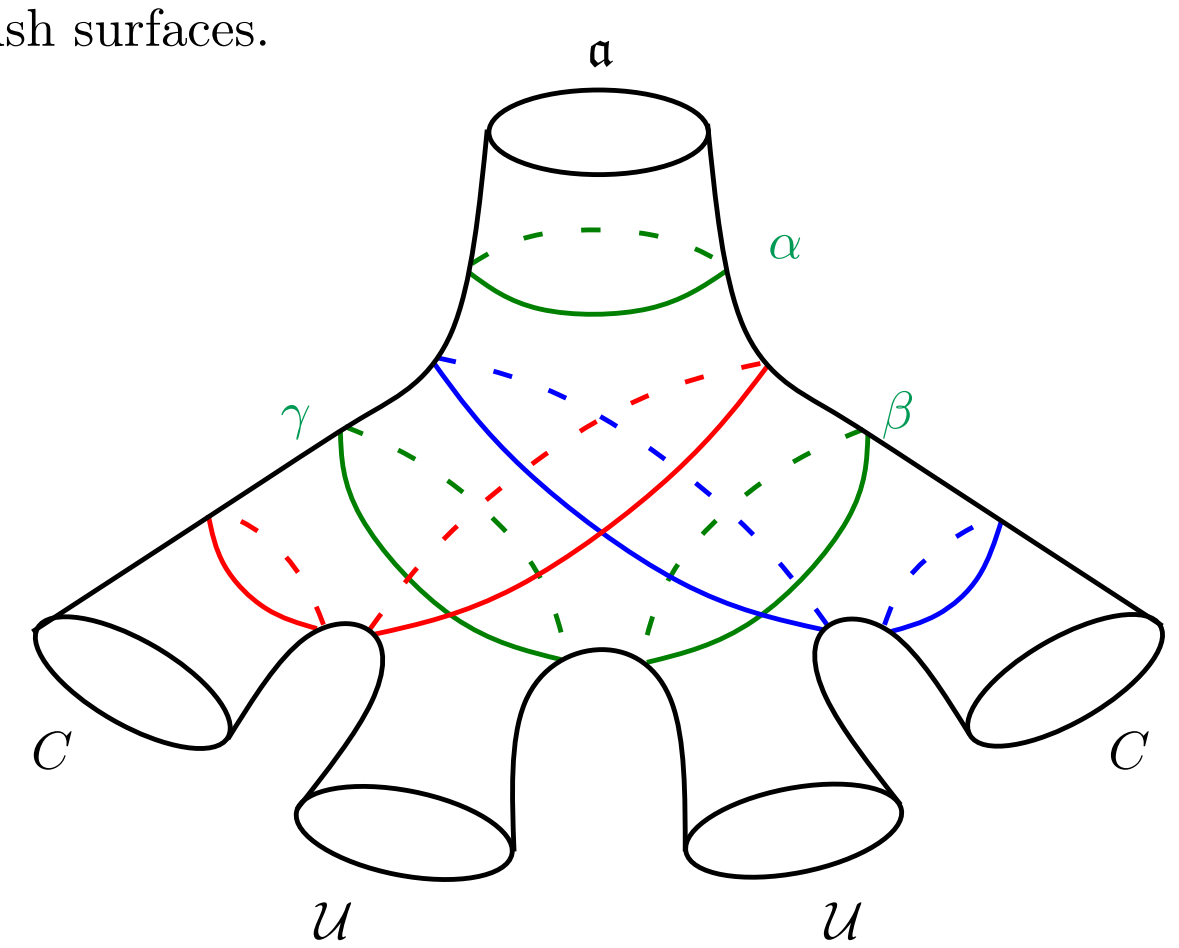
    In the case of the broom surface, recall from Subsection~\ref{subsub:endspacebroom} that $\calU$ denotes a clopen set of ends containing a copy of every non-maximal end. 
    To define the new mapping classes, we refer the reader to Figures~\ref{figure:fmaps} and \ref{figure:gmaps}, where $C$ denotes a subset of $X$ homeomorphic to a Cantor set.
    Let $f_1$ be a mapping class that, in Figure~\ref{figure:fmaps}, fixes the curve $\gamma$ and takes $\alpha$ and $\beta$ to the blue curves. Likewise, let $f_2$ be a mapping class that fixes $\beta$ and takes $\alpha$ and $\gamma$ to the red curves. Similarly, let $g_1$ be a mapping class that, in Figure~\ref{figure:gmaps}, fixes $\gamma$ and takes $\alpha$ and $\beta$ to the blue curves, and $g_2$ a mapping class that fixes $\beta$ and takes $\alpha$ and $\gamma$ to the red curves. Finally, let $h_1$ be a mapping class that, in Figure~\ref{figure:hmaps}, fixes $\alpha$ and takes $\beta$ and $\gamma$ to the blue curves, and $h_2$ a mapping class that fixes $\alpha$ and takes $\beta$ and $\gamma$ to the red curves.

    Let $P \in \calP(\Sigma)$ denote our base vertex, and consider a central pair of pants $P' \in \calP(\Sigma)$ with $d(P, P') = 1$. By definition, the boundary curves of $P$ and $P'$ intersect exactly twice, and we treat the different cases this can occur separately.

        \begin{figure}
        \centering
        \begin{overpic}[width=0.4\textwidth]{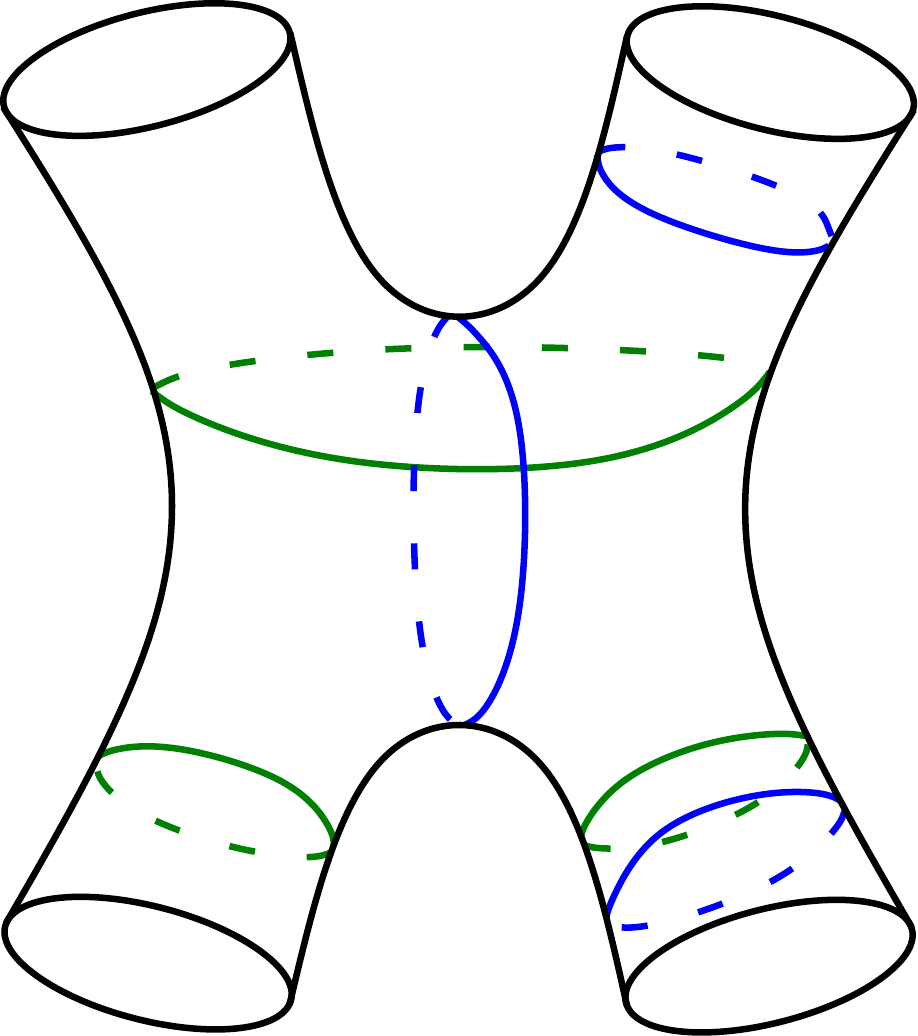}
            \put(74,60){\textcolor{ForestGreen}{$\alpha$}}
            \put(80,26.5){\textcolor{ForestGreen}{$\beta$}}
            \put(3,26.5){\textcolor{ForestGreen}{$\gamma$}}
            \put(75,100){$\Ga$}
            \put(0,100){\tiny{Unknown}}
            \put(8,-5){$C$}
            \put(75,-5){$C$}
        \end{overpic}
        \caption{When a boundary curve of $P'$ (in blue) intersects $\alpha$}
        \label{fig:int_alpha}
        \end{figure}
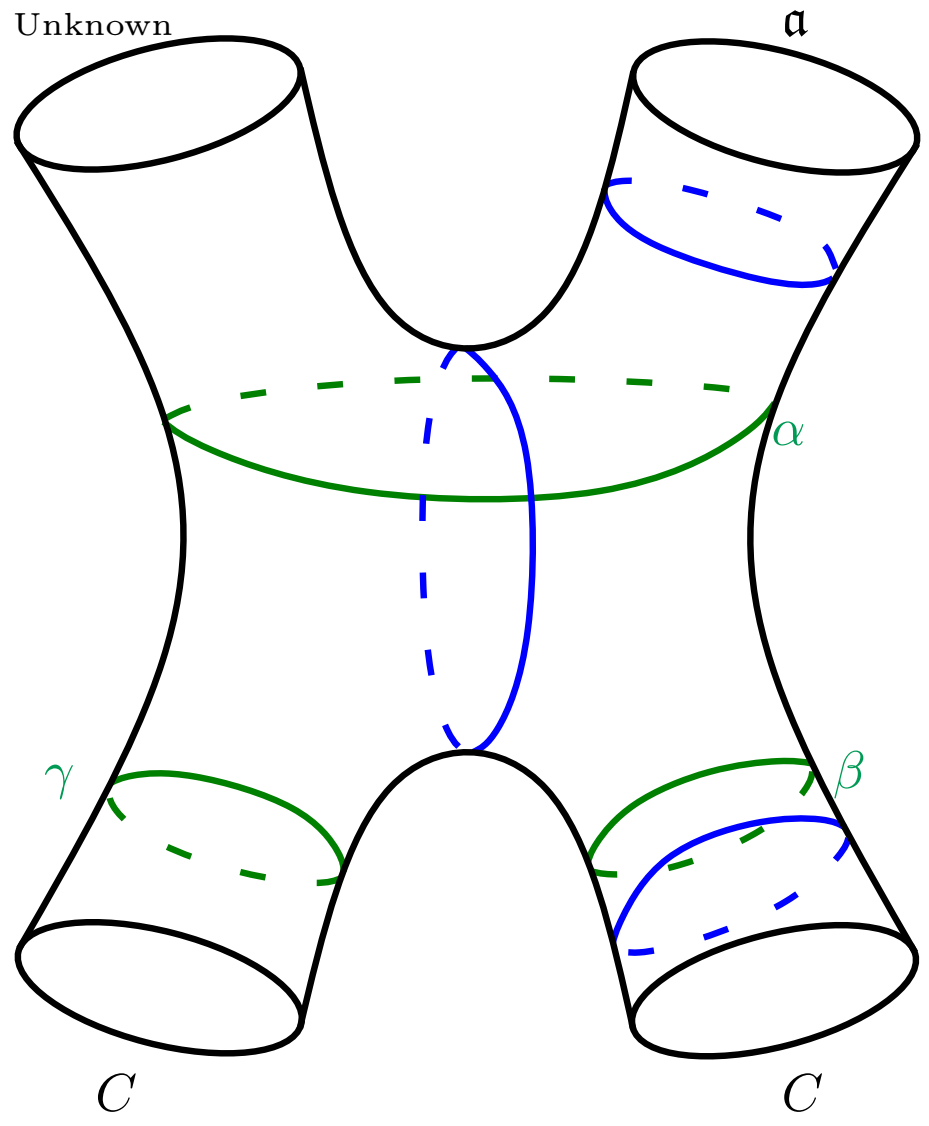

    First consider the case where one the boundary curves of $P'$ intersects the curve $\alpha$ of $P$, with the curve $\gamma$ on the same side of this curve as the singleton maximal end. This can be seen in Figure~\ref{fig:int_alpha}. Note that one complementary component of $P \cup P'$ has an arbitrary end space---we know only that it does not include a maximal end. However, the corresponding complementary component of $P \cup (f_1^{-1}P')$ has end space homeomorphic to $\calU$. Crucially, the pair of pants $g_1 P$ looks identical to $f_1^{-1}P'$, up to a homeomorphism supported on the complementary component of $P$ bounded by $\alpha$. Thus $f_1^{-1}P' \in \Stab(P) g_1 P$, and so $P' \in f_1 \Stab(P) g_1 P$. In general, any mapping class taking $P$ to $P'$ can be factored as $\phi \circ \psi$, where $\psi \in \Stab(P)$ and $\phi \in f_1 \Stab(P) g_1$, and so the mapping class taking $P$ to $P'$ is contained in $f_1 \Stab(P) g_1 \Stab(P)$.

    By an exactly analogous argument, if one of the boundary curves of $P'$ intersects the curve $\alpha$ of $P$, with the curve $\beta$ on the same side of this curve as the singleton maximal end, any mapping class taking $P$ to $P'$ is contained in $f_2 \Stab(P) g_2 \Stab(P)$.

        \begin{figure}
        \centering
        \begin{overpic}[width=0.4\textwidth]{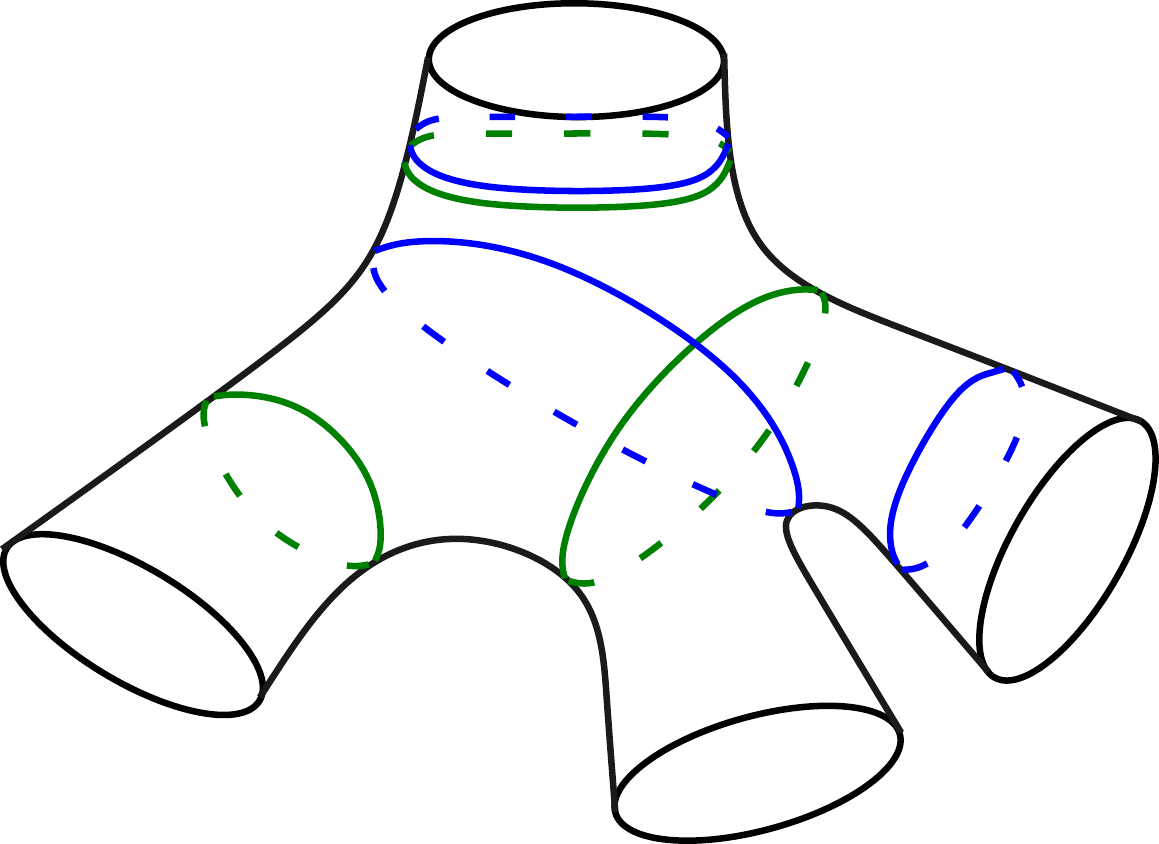}
            \put(65,59){\textcolor{ForestGreen}{$\alpha$}}
            \put(72,48.5){\textcolor{ForestGreen}{$\beta$}}
            \put(12,39){\textcolor{ForestGreen}{$\gamma$}}
            \put(48,74){$\Ga$}
            \put(55,-5){\tiny{Unknown}}
            \put(3,6){$C$}
            \put(98,18){$C$}
        \end{overpic}
        \caption{When a boundary curve of $P'$ (in blue) intersects $\beta$ and another is parallel to $\alpha$.}
        \label{fig:int_beta_par_alpha}
        \end{figure}
        
        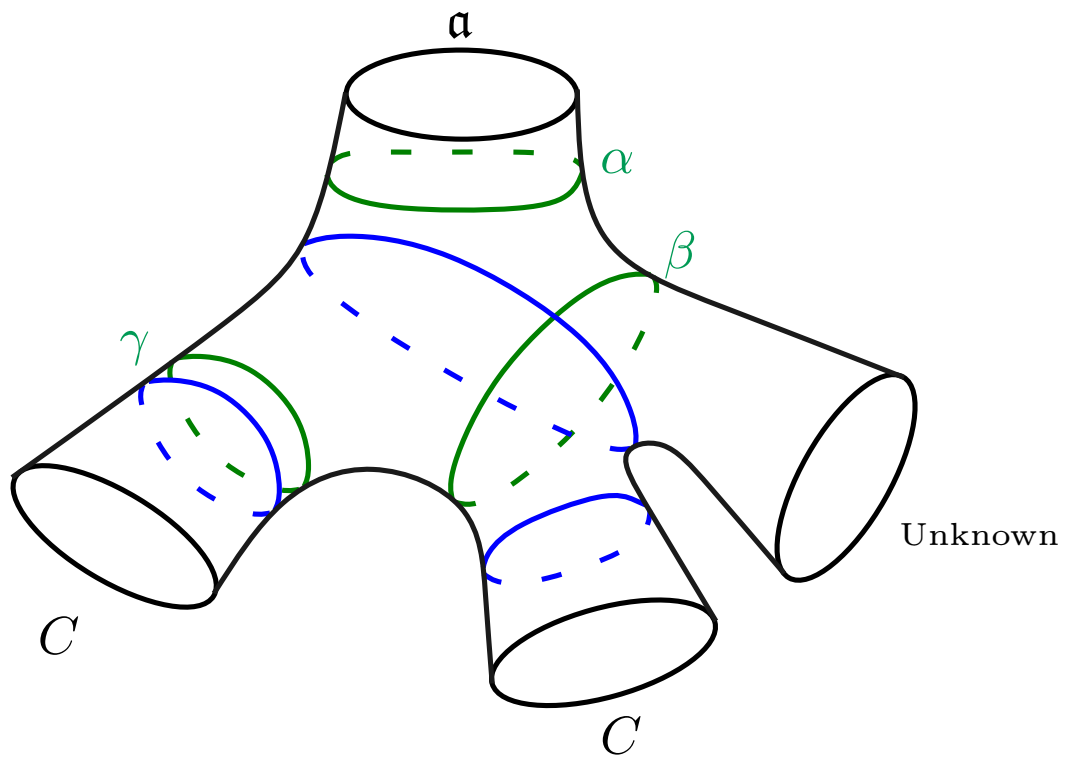
\begin{figure}
        \centering
        \begin{overpic}[width=0.4\textwidth]{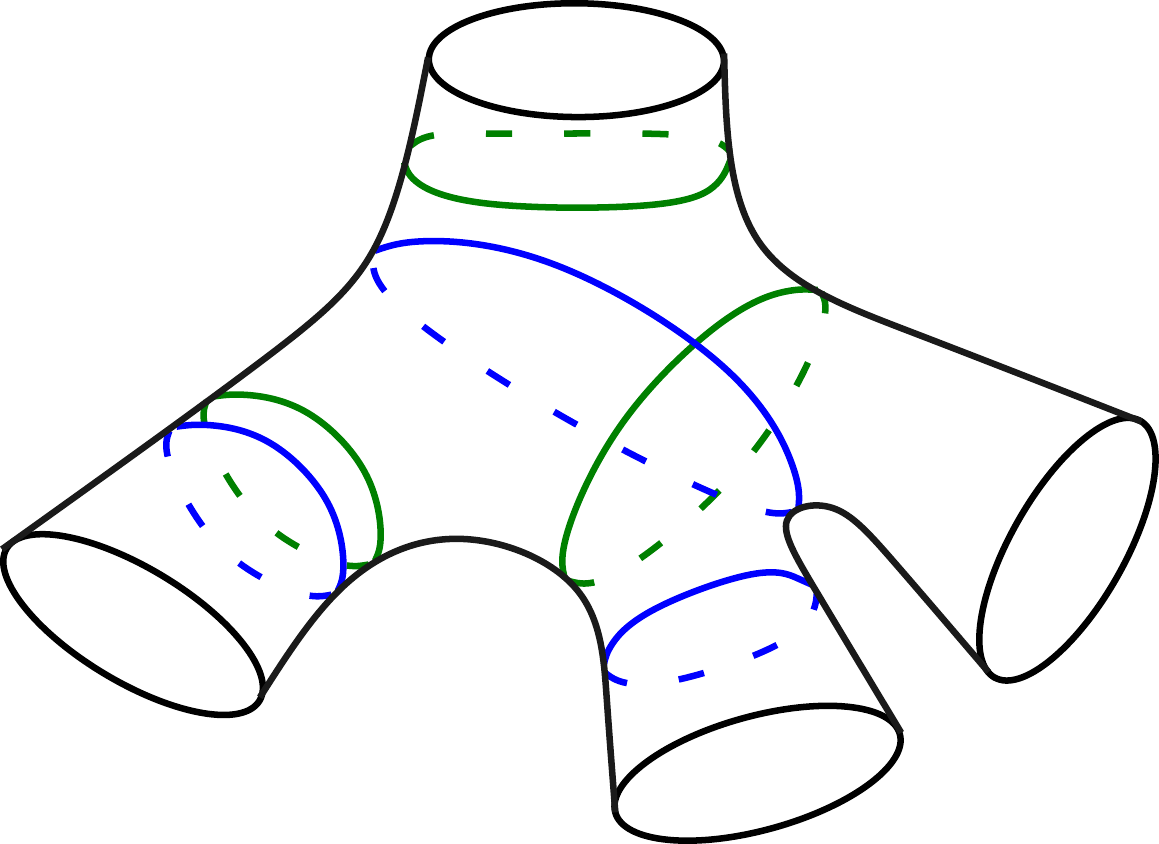}
            \put(65,59){\textcolor{ForestGreen}{$\alpha$}}
            \put(72,48.5){\textcolor{ForestGreen}{$\beta$}}
            \put(12,39){\textcolor{ForestGreen}{$\gamma$}}
            \put(48,74){$\Ga$}
            \put(65,-5){$C$}
            \put(3,6){$C$}
            \put(98,18){\tiny{Unknown}}
        \end{overpic}
        \caption{When a boundary curve of $P'$ (in blue) intersects $\beta$ and another is parallel to $\gamma$.}
        \label{fig:int_beta_par_gamma}
        \end{figure}

    Now consider the case where one of the boundary curves of $P'$ intersects the curve $\beta$ of $P$. There are two possibilities here, depicted in Figures~\ref{fig:int_beta_par_alpha} and \ref{fig:int_beta_par_gamma}, depending on which side of this curve the rest of $P'$ lies. In each of these cases we proceed as before. If $P'$ looks like Figure~\ref{fig:int_beta_par_alpha}, then $h_1^{-1}P'$ is identical to $h_2 P$ modulo a homeomorphism fixing $P$, and so any mapping class taking $P$ to $P'$ is contained in $h_1 \Stab(P) h_2 \Stab(P)$. If $P'$ looks like Figure~\ref{fig:int_beta_par_gamma}, then $g_1^{-1}P'$ is identical to $f_1 P$ modulo a homeomorphism fixing $P$, and so any mapping class taking $P$ to $P'$ is contained in $g_1 \Stab(P) f_1 \Stab(P)$.

    Finally, if one of the boundary curves of $P'$ intersects the curve $\gamma$ of $P$, by an analogous argument any mapping class taking $P$ to $P'$ is contained in $h_2 \Stab(P) h_1 \Stab(P)$ or $g_2 \Stab(P) f_2 \Stab(P)$.

    Putting this all together, let $\calF = \{\id, f_1, f_2, g_1, g_2, h_1, h_2\}$. Then the set $\{f \in \Map(\Sigma) \mid d(P, fP) \le 1\}$ is contained in $(\calF \Stab(P))^2$, which is coarsely bounded. Thus by Lemma 2.5 of \cite{SC}, the orbit map $f \mapsto fP$ is a quasi-isometry between $\Map(\Sigma)$ and $\calP(\Sigma)$.

        \begin{figure}
        \centering
        \begin{overpic}[width=0.9\textwidth]{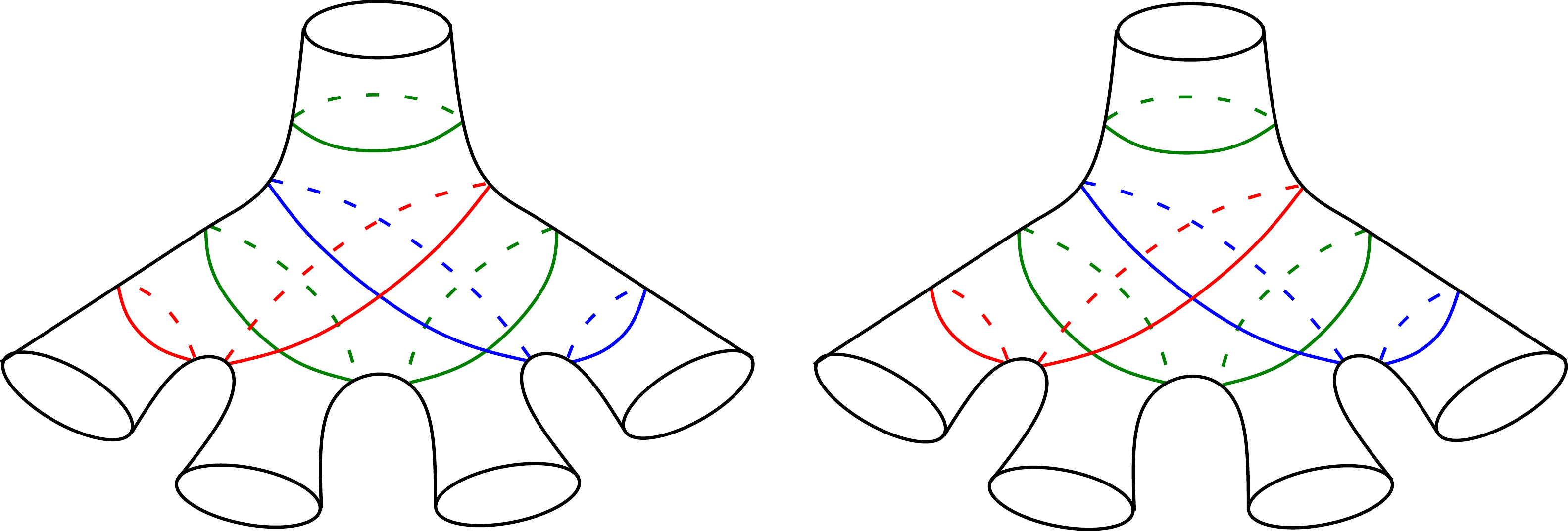}
            \put(31,26){\textcolor{ForestGreen}{$\alpha$}}
            \put(36,20){\textcolor{ForestGreen}{$\beta$}}
            \put(11,20){\textcolor{ForestGreen}{$\gamma$}}
            \put(23,34){$s$}
            \put(12.5,-2){$\calU$}
            \put(33.5,-2){$\calU$}
            \put(0.5,5){$C$}
            \put(45,5){$C$}
            \put(82.5,26){\textcolor{ForestGreen}{$\alpha$}}
            \put(87.5,20){\textcolor{ForestGreen}{$\beta$}}
            \put(62.5,20){\textcolor{ForestGreen}{$\gamma$}}
            \put(75,34){$s$}
            \put(64.5,-2){$\calV$}
            \put(84.5,-2){$\calV$}
            \put(52,5){$C$}
            \put(97,5){$C$}
        \end{overpic}
        \caption{The curves used to define the mapping classes $h_1$, $h_2$, $j_1$ and $j_2$ for the brush surface.}
        \label{figure:hjmaps}
        \end{figure}
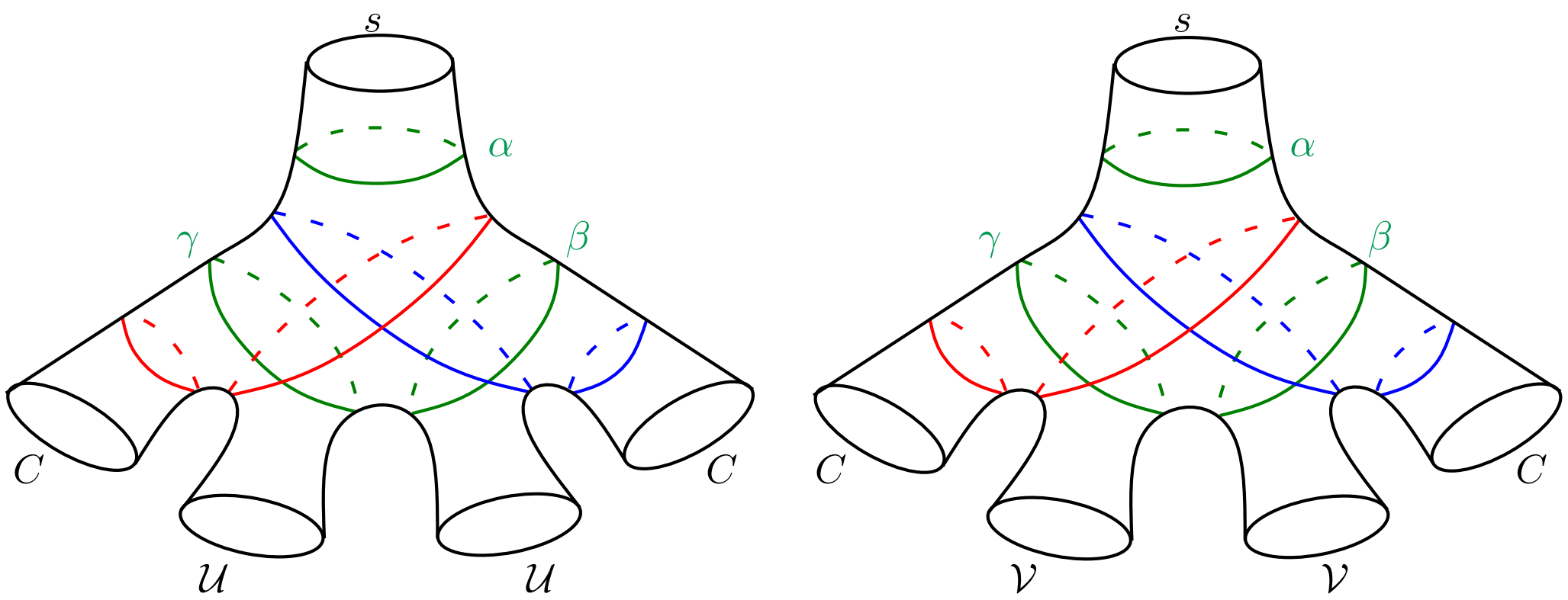
        
    In the case of the brush surface, recall that $\calU$ denotes a clopen set of ends containing a copy of every non-maximal end which is an immediate predecessor of both $\Ga$ and $\Gb$ and $\calV$ a clopen set of ends containing a copy of every non-maximal end which are immediate predecessors to the maximal end $\Gb$.
    The new mapping classes that we need to define will be based on Figures~\ref{figure:fmaps}, \ref{figure:gmaps}, and \ref{figure:hjmaps}, where $C$ will denote a subset of $E_{\Gc}$ homeomorphic to a Cantor set.
    Let $f_1$ be a mapping class that, in Figure~\ref{figure:fmaps}, fixes the curve $\gamma$ and takes $\alpha$ and $\beta$ to the blue curves. 
    Likewise, let $f_2$ be a mapping class that fixes $\beta$ and takes $\alpha$ and $\gamma$ to the red curves. 
    Similarly, let $g_1$ be a mapping class that, in Figure~\ref{figure:gmaps}, fixes $\gamma$ and takes $\alpha$ and $\beta$ to the blue curves, and $g_2$ a mapping class that fixes $\beta$ and takes $\alpha$ and $\gamma$ to the red curves.
    Let $h_1$ be a mapping class that, on the left side of Figure~\ref{figure:hjmaps}, fixes $\alpha$ and takes $\beta$ and $\gamma$ to the blue curves, and $h_2$ a mapping class that fixes $\alpha$ and takes $\beta$ and $\gamma$ to the red curves.
    Finally, let $j_1$ be a mapping class that, on the right side of Figure~\ref{figure:hjmaps}, fixes $\alpha$ and takes $\beta$ and $\gamma$ to the blue curves, and $j_2$ a mapping class that fixes $\alpha$ and takes $\beta$ and $\gamma$ to the red curves.

    The case where one of the boundary curves of $P' \in \calP(\Sigma)$ intersects the curve $\alpha$ of $P$ follows exactly as in the case of the broom.
    We now consider the case where one of the boundary curves of $P'$ intersects the curve $\beta$ of $P$. There are two possibilities here, depicted in Figures~\ref{fig:int_beta_par_alpha} and \ref{fig:int_beta_par_gamma}, depending on which side of this curve the rest of $P'$ lies. In each of these cases we proceed as before. If $P'$ looks like Figure~\ref{fig:int_beta_par_alpha}, then either $h_1^{-1}P'$ is identical to $h_2 P$ modulo a homeomorphism fixing $P$, or $j_1^{-1}P'$ is identical to $j_2 P$.
    Therefore, any mapping class taking $P$ to $P'$ is contained in $h_1 \Stab(P) h_2 \Stab(P)$ or in $j_1 \Stab(P) j_2 \Stab(P)$. If $P'$ looks like Figure~\ref{fig:int_beta_par_gamma}, then $g_1^{-1}P'$ is identical to $f_1 P$ modulo a homeomorphism fixing $P$, and so any mapping class taking $P$ to $P'$ is contained in $g_1 \Stab(P) f_1 \Stab(P)$.

    Finally, if one of the boundary curves of $P'$ intersects the curve $\gamma$ of $P$, an analogous argument shows that any mapping class taking $P$ to $P'$ is contained in $h_2 \Stab(P) h_1 \Stab(P)$, $j_2 \Stab(P) j_1 \Stab(P)$, or $g_2 \Stab(P) f_2 \Stab(P)$.

    Putting this all together, let $\calF = \{\id, f_1, f_2, g_1, g_2, h_1, h_2, j_1, j_2\}$. Then the set $\{f \in \Map(\Sigma) \mid d(P, fP) \le 1\}$ is contained in $(\calF \Stab(P))^2$, which is coarsely bounded. Thus by Lemma 2.5 of \cite{SC}, the orbit map $f \mapsto fP$ is a quasi-isometry between $\Map(\Sigma)$ and $\calP(\Sigma)$.

    Finally, we analyze the case of the tripod.
    Recall from Subsection~\ref{subsub:endspacetripod} that $\calU$ denotes a clopen set of ends containing a copy of every non-maximal end.
    The new mapping classes that we need to define will be based on Figure~\ref{figure:fmapstripod}.
    Let $f_1$ be a mapping class that, in Figure~\ref{figure:fmapstripod}, fixes the curve $\gamma$ and takes $\alpha$ and $\beta$ to the blue curves. Likewise, let $f_2$ be a mapping class that fixes $\beta$ and takes $\alpha$ and $\gamma$ to the red curves.
    By changing the roles of the labels of the maximal ends, we may also construct the maps $g_1, g_2, h_1,$ and $h_2$.
    \begin{figure}
        \centering
        \begin{overpic}[width=0.5\textwidth]{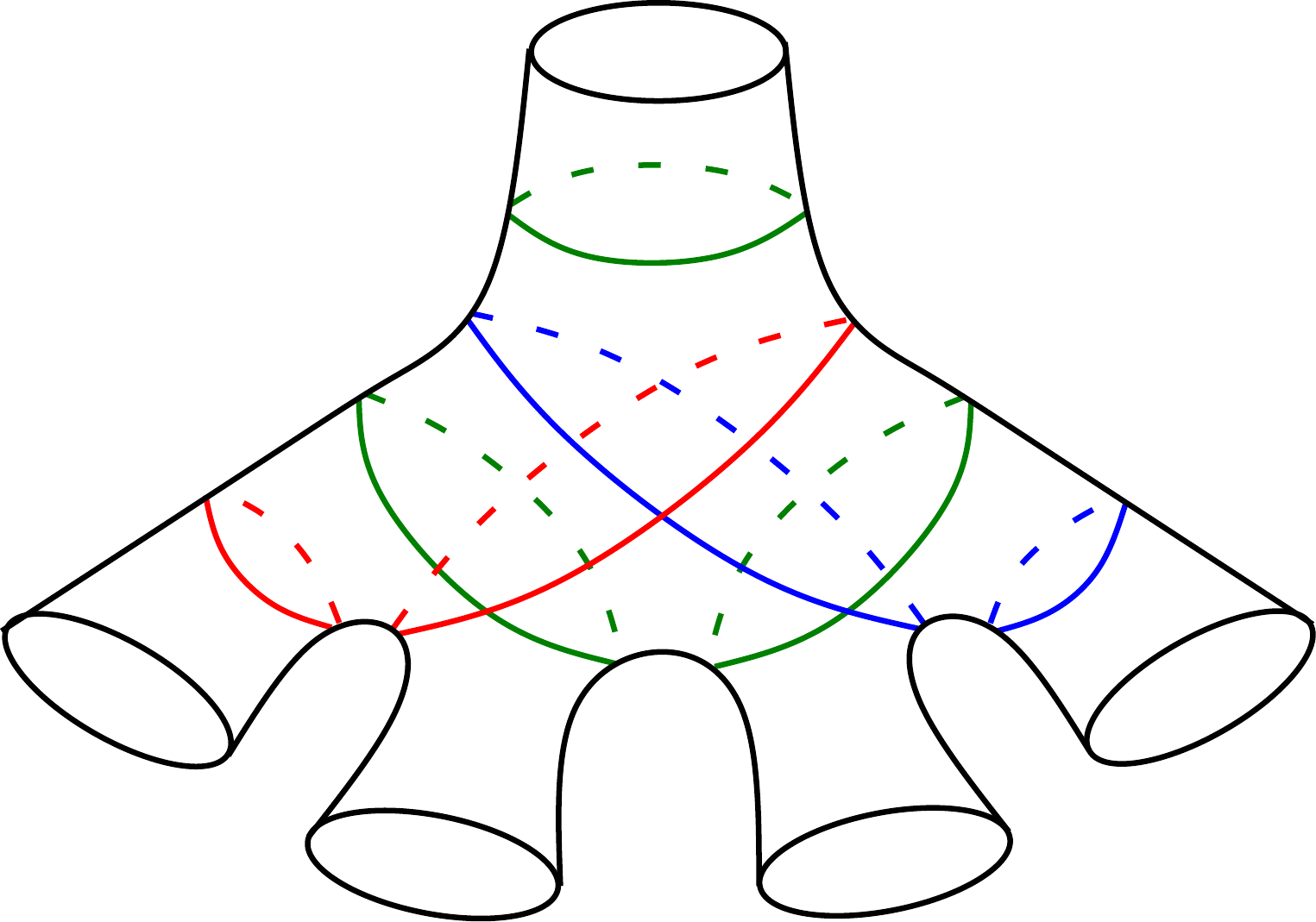}
            \put(65,55){\textcolor{ForestGreen}{$\alpha$}}
            \put(75,40){\textcolor{ForestGreen}{$\beta$}}
            \put(22,40){\textcolor{ForestGreen}{$\gamma$}}
            \put(49,72){$\Ga$}
            \put(25,-5){$\calU$}
            \put(70,-5){$\calU$}
            \put(0,10){$\Gc$}
            \put(95,10){$\Gb$}
        \end{overpic}
        \caption{The curves used to define the mapping classes $f_1$ and $f_2$ for the tripod surface.}
        \label{figure:fmapstripod}
        \end{figure}

    Consider a central pair of pants $P'$ with $d(P, P') = 1$. By definition, the boundary curves of $P$ and $P'$ intersect exactly twice. We cover each of the cases where the boundary curves of $P$ and $P'$ intersect exactly twice individually.

        \begin{figure}
        \centering
        \begin{overpic}[width=0.4\textwidth]{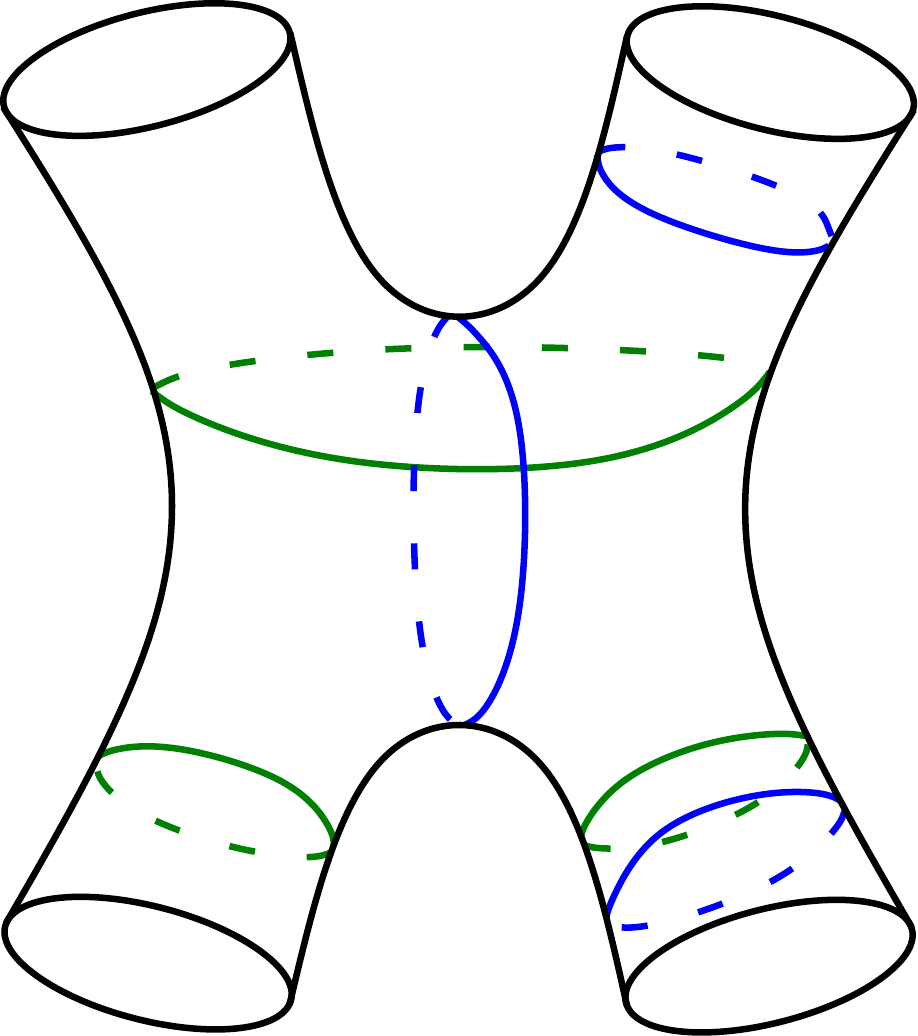}
            \put(74,60){\textcolor{ForestGreen}{$\alpha$}}
            \put(80,26.5){\textcolor{ForestGreen}{$\beta$}}
            \put(3,26.5){\textcolor{ForestGreen}{$\gamma$}}
            \put(75,100){$\Ga$}
            \put(0,100){\tiny{Unknown}}
            \put(8,-5){$\Gb$}
            \put(75,-5){$\Gc$}
        \end{overpic}
        \caption{When a boundary curve of $y$ (in blue) intersects $\alpha$}
        \label{fig:int_alphaTripod}
        \end{figure}

    First consider the case where one the boundary curves of $P'$ intersects the curve $\alpha$ of $P$, with the curve $\gamma$ on the same side of this curve as the singleton maximal end. This can be seen in Figure~\ref{fig:int_alphaTripod}. Note that one complementary component of $P \cup P'$ has an arbitrary end space---we know only that it does not include a maximal end. However, the corresponding complementary component of $P \cup (f_1^{-1}P')$ has end space homeomorphic to $\calU$. Crucially, the pair of pants $g_1 P$ looks identical to $f_1^{-1}P'$, up to a homeomorphism supported on the complementary component of $P$ bounded by $\alpha$. Thus $f_1^{-1}P' \in \Stab(P) g_1 P$, and so $P' \in f_1 \Stab(P) g_1 P$. In general, any mapping class taking $P$ to $P'$ can be factored as $\phi \circ \psi$, where $\psi \in \Stab(P)$ and $\phi \in f_1 \Stab(P) g_1$, and so the mapping class taking $P$ to $P'$ is contained in $f_1 \Stab(P) g_1 \Stab(P)$.
    The remaining cases are done analogously, but changing the labeling in Figure~\ref{fig:int_alphaTripod}.
    
    Putting this all together, let $\calF = \{\id, f_1, f_2, g_1, g_2, h_1, h_2\}$. Then the set $\{f \in \Map(\Sigma) \mid d(P, fP) \le 1\}$ is contained in $(\calF \Stab(P))^2$, which is coarsely bounded. Thus by Lemma 2.5 of \cite{SC}, the orbit map $f \mapsto fP$ is a quasi-isometry between $\Map(\Sigma)$ and $\calP(\Sigma)$.

    Thus, we have shown that $\Map(\Sigma)$ is quasi-isometric to $\calP(\Sigma)$ where $\Sigma$ is either the broom, brush, or tripod surfaces.
\end{proof}

Applying Proposition~\ref{prop:qi}, Corollary~\ref{Cor:HypP}, and the work of Gromov which shows hyperbolicity is invariant under quasi-isometries \cite{Gromov1987}, we have proved that the mapping class group of the broom, brush and tripod surfaces is hyperbolic.

\begin{proposition}\label{Prop:MapHyp}
    Suppose that $\Sigma$ is either the broom, brush or tripod surface.
    Then $\Map(\Sigma)$ is non-trivially $\delta$-hyperbolic, where the hyperbolicity constant is independent of the underlying surface.
\end{proposition}

\section{Infinite Coarse Rank}

In this section, we will complete our analysis of the surfaces where the end complexity, $\zeta(\Sigma) = 3$.
The main tool in this section is to show that for any $k \in \NN$, there exists $k$ pairwise disjoint non-peripheral curves in each of the remaining surfaces with $\zeta(\Sigma) = 3$. 
Recall that a subsurface $S\subset \Sigma$ is called \textit{peripheral} if for any compact $K \subset \Sigma$, there exists 
$f \in \Map(\Sigma)$ such that $K \cap f(S)=\emptyset$. Otherwise, a surface is \textit{nonperipheral}. 
Equivalently, a surface is peripheral if there exists some end $x$ of $\Sigma$ (or puncture) such that for any open $x \in U\subset \Sigma$, we can find $f\in \Map(\Sigma)$ such that $f(S) \subset U$. Morally, a surface is peripheral if it can be pushed off to an end.
Once we have found the required disjoint nonperipheral curves, we apply Proposition 3.16 of \cite{BNQR} which implies that $\Map(\Sigma)$ contains a quasi-isometrically embedded copy of $\ZZ^k$, which proves our claim.

In this section, we will proceed by first analyzing the surfaces for which $\zeta = 3$ that are not broom-, brush-, or tripod-like, then analyze the surfaces which are broom-, brush-, or tripod-like with finite genus, and ending with the surfaces which are broom-, brush-, or tripod-like with infinite genus.

\subsection{Surfaces with $\zeta = 3$ that are not broom-, brush-, or tripod-like}
For brevity, we define a new invariant.
\begin{definition}
Suppose that $\Sigma$ is a CB-generated 
infinite-type surface. Define:
$\xi(\Sigma)$ = 2 (number of Cantor-type maximal ends)
 + (number of isolated maximal ends)
 + (number of immediate predecessors that do not accumulate to every maximal end).
 In particular, $\xi(\Sigma) \ge \zeta(\Sigma)$, and 
 $\xi(\Sigma) = 3$ if and only if $\Sigma$ is broom-like, brush-like, or tripod-like.
 \end{definition}

\begin{lemma}
    Suppose that $\xi(\Sigma) > 3$. 
    Then $\Map(\Sigma)$ has infinite coarse rank.
\end{lemma}
\begin{proof}
We restrict our attention to surfaces with no 
essential shift and $\zeta(\Sigma) = 3$, as the remaining cases are treated in \cite{BNQR}. 
We treat different cases of end spaces separately. 
\begin{itemize}
\item If $\calM(\Sigma) = \{\Ga, \Gb, \Ge\}$ then without loss of generality,
there must exist some predecessor end with equivalence class $E$ that does not
converge to $\Ga$. If $E$ were countable, then $\Sigma$ 
is not CB generated, so $E$ must be Cantor-type. Let $E = E_1\sqcup E_2$ be 
a clopen decomposition of $E$ such that $E_1$ does not 
converge to any maximal end. Let 
$\alpha_0$ be a simple closed 
curve separating $\Ga \cup E_1$ from $E_2 \cup \Gb \cup \Ge$. Then $\alpha_0$ is 
non-peripheral (for example, any orbit of $\alpha_0$ will intersect a pair of pants separating the maximal ends). Let $S_{\alpha_0}$ be the surface cut out by $\alpha_0$ which contains $\Ga$. Notice now that there exists an $\alpha_1 \subset S_{\alpha_0}$ with the same properties as $\alpha_0$ and hence is non-peripheral. Proceeding in this manner 
constructs infinitely-many disjoint non-peripheral curves.
\item If $\calM(\Sigma) = \{\Ga, \Gb\}$ then 
let $E_{\Ge}$ be the equivalence class of ends that is an immediate predecessor of $\Ga$ but not $\Gb$, and note that there must exist another Cantor-type equivalence class, $E$ 
which converges to one maximal end but not the other. 
The construction of infinitely-many non-peripheral curves follows analogously, this time treating $E_{\Ge}$ in the same way as we treated $\Ge$ in the previous case.
\item If $\calM(\Sigma) = X \cup \Ga$, where $X$ is a 
Cantor set, then as before, there exists a Cantor-type equivalence class $E$ such that $E$ does not converge to either $X$ or $E$ does not converge to $\Ga$. 
Suppose first that $E$ does not converge to $X$.
We may then write $E = E_1\sqcup E_2$ and $X = X_1\sqcup X_2$ (also a clopen decomposition) and define $\alpha_0$ 
to be a simple closed curve separating $E_1 \cup X_1$ from $E_2 \cup X_2 \cup \Ga$. This curve is non-peripheral, as it intersects all four-times punctured spheres which separate $X$ and $E$, but it does not have a complementary component with ends from both of these sets since this property is mapping class invariant (ie. any orbit of this subsurface will intersect $\alpha_0$). 
Notice now that there exists an $\alpha_1$ in the subsurface of $\Sigma \ssm \alpha_0$ containing $\Ga$ with the same properties as $\alpha_0$. 
Proceeding in this manner constructs infinitely-many disjoint non-peripheral curves. 

Now suppose that $E$ converges to $X$ but not to $\Ga$.
Define $\alpha_0$ to be a simple closed curve separating $\Ga \cup E_1$ from $E_2 \cup X$. This 
curve intersects all pairs of pants $P$ which contain $\Ga$ but not $E$ in a complementary component of $P$, and parts of $X$ in all other complementary components. Applying a similar argument to before, we can construct infinitely-many such curves.
\end{itemize}

Thus, when $\xi >3$, $\Sigma$ admits infinitely-many disjoint non-peripheral surfaces.
\end{proof}

\subsection{Broom-like, brush-like, and tripod-like surfaces with finite genus}
In this subsection, we will show that broom-like, brush-like, and tripod-like surfaces with finite genus have infinite coarse rank.
To this end, we will require an \textit{anchor surface} of $\Sigma$, as introduced by Bar-Natan--Qing--Rafi, which is a subsurface whose boundary components define the partition of the endspace as introduced by Mann--Rafi (see Proposition 5.5 of \cite{MR}).
An anchor surface $K_0$ is said to be \textit{minimal} if the genus and number of components of $\Sigma \ssm K_0$ are both minimal. Using the anchor surface, we produce a collection of $k$ pairwise disjoint non-peripheral curves and finish by invoking Proposition 3.16 of \cite{BNQR}.

\begin{lemma}\label{Lem:InfCoarseRankFinteGenus}
    Let $\Sigma$ be a stable, infinite-type surface with finite positive genus. Suppose further that $\Map(\Sigma)$ is CB-generated.
    Then $\Map(\Sigma)$ has infinite coarse rank.
\end{lemma}

\begin{proof}
    By definition of a minimal anchor surface, since the surface $\Sigma$ has finite genus, the genus of the anchor surface must be equal to the genus of the surface.
    Let $\gamma$ denote a simple closed curve which separates the finite genus from the rest of the surface. 
    Notice that the subsurface of $\Sigma \ssm \gamma$ containing the finite genus, denoted $S_{\gamma,g}$ is a non-peripheral subsurface.
    Indeed, suppose that there exists some $f \in \Map(\Sigma)$ such that $f(S_{\gamma,g}) \cap K = \emptyset$.
    But then $S_{\gamma,g}$ must be mapped to a subsurface of $\Sigma \ssm K$ which contains genus, a contradiction since $S_{\gamma,g}$ contains all the genus of the surface.

    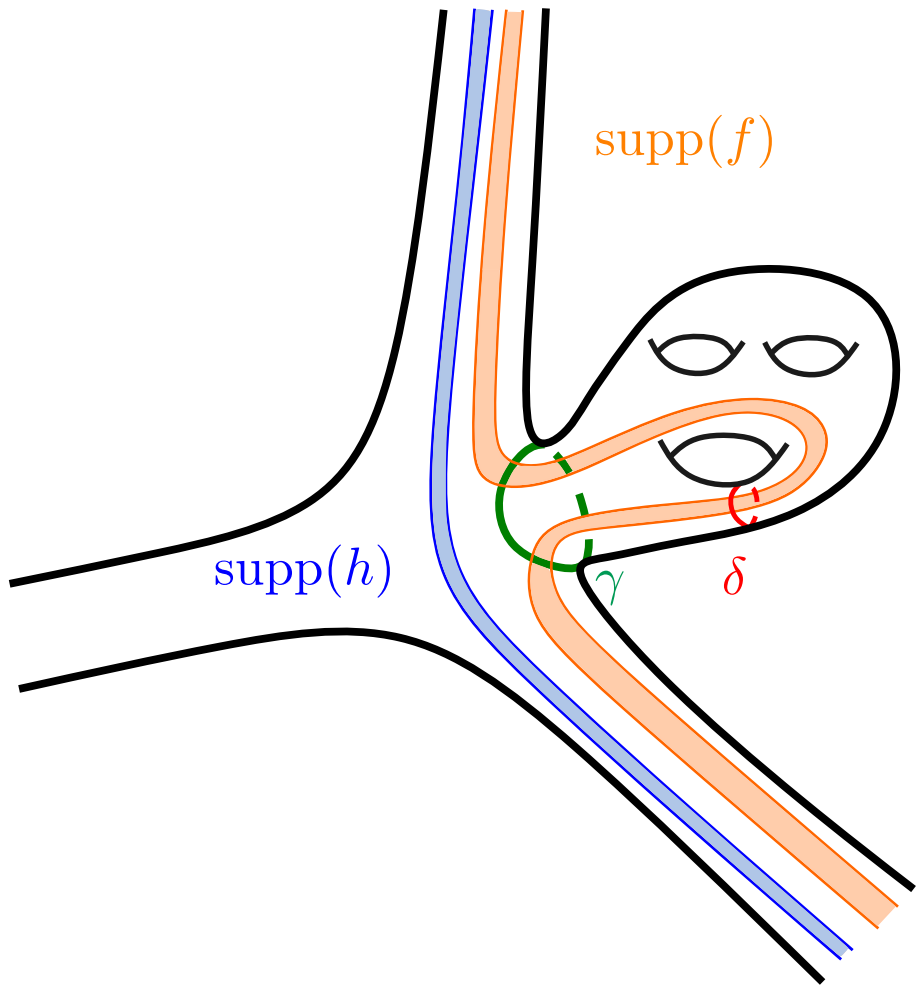
\begin{figure}
        \centering
        \begin{overpic}[width=0.4\textwidth]{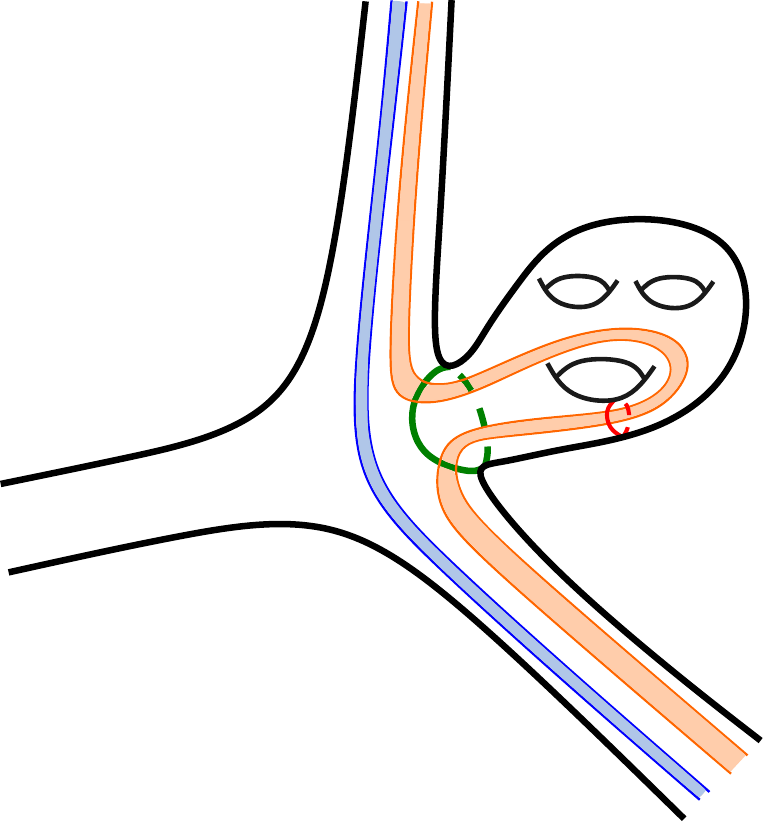}
            \put(60,40){\textcolor{ForestGreen}{$\gamma$}}
            \put(73,40){\textcolor{red}{$\delta$}}
            \put(60,85){\textcolor{orange}{$\mathrm{supp}(f)$}}
            \put(21,41){\textcolor{blue}{$\mathrm{supp}(h)$}}
        \end{overpic}
        \caption{Modifying $\mathrm{supp}(f)$ so that it intersect $S_{\gamma,g}$ and $\delta$.}
        \label{fig:shiftnonperipheralcurves}
        \end{figure}

    We now note that Section 6 of \cite{MR}, that there exist non-essential shifts in $f \in \Sigma$.
    Should $S_{\gamma,g} \cap \mathrm{supp}(f) = \emptyset$, then we modify $\mathrm{supp}(f)$ to non-trivially intersect $S_{\gamma,g}$ to construct a shift $h \in \Map(\Sigma)$.
    See Figure~\ref{fig:shiftnonperipheralcurves}.
    We now select a non-separating simple closed curve $\delta$ in $S_{\gamma,g}$ which non-trivially intersects $\mathrm{supp}(f)$, which exists due to the existence of at least one genus in $S_{\gamma,g}$.
    Fix $j \in \NN$.
    Notice that $h^i(\delta)$ for $1 \leq i \leq k$ forms a collection of $k$ pairwise disjoint nonperipheral curves in $\Sigma$.
    Thus, by Proposition 3.16 of \cite{BNQR}, $\Map(\Sigma)$ has infinite coarse rank.
\end{proof}

\subsection{Broom-like, brush-like, and tripod-like surfaces with infinite genus}

In this subsection, we prove that broom-like, brush-like, and tripod-like surfaces with infinite genus have infinite coarse rank.

\begin{lemma}\label{Lem:InfCoarseRankInfinteGenus}
    Let $\Sigma$ be either a tripod-like, broom-like, or brush-like surface with infinite genus, where at least one maximal end is not accumulated by genus.
    Then $\Map(\Sigma)$ has infinite coarse rank.
\end{lemma}

\begin{proof}

    Let $K_0$ denote the minimal anchor surface for $\Sigma$.
    Notice that every complementary component of $K_0$ either has finite genus, infinite genus, or no genus.
    Let $z \in \calM(\Sigma)$ denote one of the maximal ends which has a neighborhood not containing infinite genus.
    Let $k \in \NN$.
    Let $\gamma$ be a separating curve such that $\Sigma \ssm \gamma$ contains a subsurface which contains $z$ and has genus equal to $j$.
    Let this subsurface of $\Sigma \ssm \gamma$ be denoted by $S_{\gamma,j}$.
    We will show that $\gamma$ is a non-peripheral curve.
    Let $K$ denote a compact connected surface which contains both the minimal anchor surface, $K_0$, and the finite genus contained in the surface $S_{\gamma,j}$.
    As with $K_0$, the complementary components of $K$ either have infinite genus or no genus.

    Suppose that there exists from $f \in \Map(\Sigma)$ such that $f(\gamma) \cap K = \emptyset$.
    Therefore, by definition of $K$, $\gamma$ must be mapped to a complementary surface of $\Sigma \ssm K$ containing infinite genus.
    This implies that $z$ is mapped into a component of $\Sigma \ssm K$ containing infinite genus, which is a contradiction as then the image of $K_0$ under $f$ is unable to partition the endspace as by the definition of an anchor surface.
    Thus $\gamma$ is a non-peripheral curve. By this argument, one can create a ``stack" of $j$ disjoint separating curves for any $j \in \NN$.
    See Figure~\ref{fig:infinitegenuscoarserank} for an illustration of $j$ disjoint separating curves in the case of the tripod-like surface.
        \begin{figure}
        \centering
        \begin{overpic}[width=0.25\textwidth]{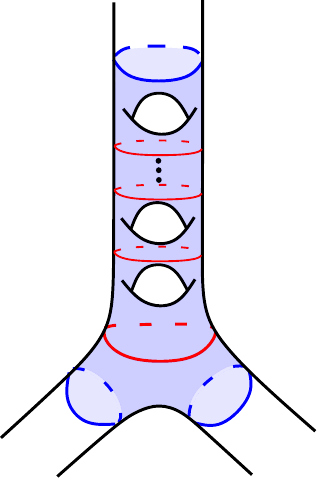}
            \put(44,68){\textcolor{red}{$\gamma_1$}}
            \put(44,57){\textcolor{red}{$\gamma_{j-2}$}}
            \put(44,45){\textcolor{red}{$\gamma_{j-1}$}}
            \put(46,31){\textcolor{red}{$\gamma_j$}}
            \put(17,78){\textcolor{blue}{$K$}}
            \put(32,100){$x$}
        \end{overpic}
        \caption{$\Sigma$ is a tripod-like surface with infinite genus, $K \subset \Sigma$ a compact subsurface containing $j$ disjoint non-peripheral separating curves.}
        \label{fig:infinitegenuscoarserank}
        \end{figure}
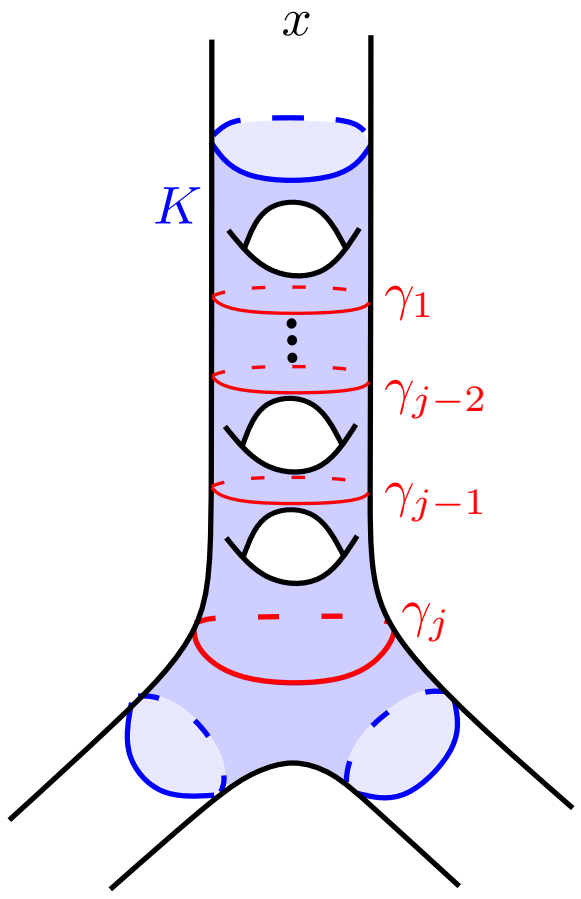
    Thus, by Proposition 3.16 of \cite{BNQR}, $\Map(\Sigma)$ has infinite coarse rank.
\end{proof}

Applying Proposition~\ref{Prop:MapHyp}, Lemmas~\ref{Lem:InfCoarseRankFinteGenus} and \ref{Lem:InfCoarseRankInfinteGenus}, together with the work of Bar-Natan--Qing--Rafi, we have proven the main theorem of this paper:\\

\noindent \textbf{Theorem 1.1. }
\textit{Suppose that $\Sigma$ is stable and let $\Map(\Sigma)$ have a coarsely bounded generating set. Then $\Map(\Sigma)$ has infinite coarse rank or is hyperbolic. Moreover, If $\Map(\Sigma)$ is hyperbolic and is not coarsely bounded, then $\Sigma$ is either a broom, brush, or tripod surface.}

\bibliographystyle{plain}
\bibliography{bib}

@article{Bav ,
author = {Bavard, Juliette},
title = {Hyperbolicit\'{e} du graphe des rayons et quasi-morphismes sur un gros groupe modulaire},
journal = {Geometry \& Topology},
volume = {20},
year = {2016},
number = {1},
pages = {491-535}
}

@article{Bow,
author = {Bowditch, Brian},
year = {2014},
month = {07},
pages = {},
title = {Uniform hyperbolicity of the curve graphs},
volume = {269},
journal = {Pacific Journal of Mathematics},
doi = {10.2140/pjm.2014.269.269}
}

@article{Bow2,
author = {Bowditch, Brian},
year = {2006},
title = {Intersection numbers and the hyperbolicity of the curve complex},
volume = {2006},
number={598},
pages={105-129},
journal = {Journal f\"{u}r die reine und angewandte Mathematik},
doi = {10.2140/pjm.2014.269.269}
}

@misc{GRV,
author = {Grant, Curtis and Rafi, Kasra and Verberne, Yvon},
title = {Asymptotic dimension of big mapping class groups},
year={2021},
note={arXiv:2110.03087}}

@misc{MR ,
author = {Mann, Katie and Rafi, Kasra},
title = {Large scale geometry of big mapping class groups},
journal = {arXiv:1912.10914},
year = {2019},
pages = {1-41}}

@article{R, 
title={On the Classification of Noncompact Surfaces},
volume={106}, DOI={10.1090/S0002-9947-1963-0143186-0}, 
journal={Trans. Amer. Math. Soc.}, 
author={Ian Richards}, 
year={1963}, 
pages={259-269}
}

@misc{Rose, title={{Coarse Geometry of Topological Groups. Book Manuscript}}, url={http://homepages.math.uic.edu/~rosendal/PapersWebsite/Coarse-Geometry-Book17.pdf}, author={Rosendal, Christian}}

@misc{BNQR,
author = {Bar-Natan Assaf and Qing, Yulan and Rafi, Kasra},
title = {The non-peripheral curve graph and divergence in big mapping class groups},
journal = {arXiv:2603.21560},
year = {2026},
pages = {1-37}}

@misc{SC,
  doi = {10.48550/ARXIV.2006.14760},  
  url = {https://arxiv.org/abs/2006.14760},  
  author = {Schaffer-Cohen, Anschel},  
  title = {Graphs of curves and arcs quasi-isometric to big mapping class groups},  
  publisher = {arXiv},  
  year = {2020},  
  copyright = {arXiv.org perpetual, non-exclusive license}
}

@article{HPW,
author = {Hensel, Sebastian and Przytycki, Piotr and Webb, Richard},
year = {2013},
month = {01},
pages = {},
title = {Slim unicorns and uniform hyperbolicity for arc graphs and curve graphs}
}

@book{BKR,
  editor    = {Baykur, R. {\.I}nan{\c{c}} and Kirby, Robion C. and Ruberman, Daniel},
  title     = {{K3: A New Problem List in Low-Dimensional Topology}},
  series    = {Mathematical Surveys and Monographs},
  volume    = {295},
  publisher = {American Mathematical Society},
  address   = {Providence, RI},
  year      = {2026},
  pages     = {xi+430},
}

@incollection{Gromov1987,
  author    = {Gromov, Mikhail},
  title     = {Hyperbolic groups},
  booktitle = {Essays in group theory},
  series    = {Mathematical Sciences Research Institute Publications},
  volume    = {8},
  editor    = {Gersten, Stephen M.},
  pages     = {75--263},
  publisher = {Springer},
  address   = {New York},
  year      = {1987}
}

@misc{KS,
  title         = {Geometric models and asymptotic dimension for infinite-type surface mapping class groups},
  author        = {Kopreski, Michael C. and Shaji, George},
  year          = {2025},
  eprint        = {2508.06679},
  archivePrefix = {arXiv},
  primaryClass  = {math.GT},
  doi           = {10.48550/arXiv.2508.06679},
  url           = {https://arxiv.org/abs/2508.06679}
}

\end{document}